\documentclass[leqno, dvips, 12pt]{amsart}
\usepackage{amssymb}
\usepackage{subfigure}
\usepackage{graphics}
\usepackage{float}
\usepackage{color}
\newcommand{\mbf}[1]{\ensuremath{\mathbf{#1}}}
\newcommand{\mbb}[1]{\ensuremath{\mathbb{#1}}}
\usepackage{euscript}
\usepackage{verbatim}
\newcommand{\mop}{\operatorname}
\newcommand{\E}[1]{\ensuremath{\EuScript{#1}}}
\renewcommand{\hom}[2]{\ensuremath{\mop{\textnormal{Hom}}(#1,#2)}}
\newcommand{\hpg}{\hom{\pi}{G}}
\newcommand{\hpgmg}{\ensuremath{\hpg/G}}
\newcommand{\lie}[2]{\ensuremath{\mop{#1}(2, \mathbb{#2})}}
\newcommand{\slr}{\lie{SL}{R}}
\newcommand{\slc}{\lie{SL}{C}}
\newcommand{\pglr}{\lie{PGL}{R}}
\newcommand{\glr}{\lie{GL}{R}}
\newcommand{\pglz}{\lie{PGL}{Z}}
\newcommand{\glz}{\lie{GL}{Z}}
\newcommand{\MCG}{{\pi_0(\mop{Homeo}(M))}}
\newcommand{\Z}{{\mathbb Z}}
\newcommand{\C}{{\mathbb C}}
\newcommand{\R}{{\mathbb R}}
\newcommand{\rep}[2]{\ensuremath{\E{R}_{#1,#2}}}
\newcommand{\cv}[2]{\ensuremath{\E{X}_{#1,#2}}}
\newcommand{\PMB}{\ensuremath{S_{0,2}^{\sharp}}}
\newcommand{\PKB}{\ensuremath{S_{1,1}^{\sharp}}}
\newcommand{\pants}{\ensuremath{S_{0,3}}}

\newcommand{\qps}{\ensuremath{S_{0,4}}}
\newcommand{\PMBt}{\ensuremath{\widetilde{S_{0,2}^{\sharp}}}}
\newcommand{\inv}{^{-1}}
\newcommand{\oM}{\ensuremath{\Omega_0^M}}
\newcommand{\ec}{\ensuremath{\E{E}_c}}
\newcommand{\zbar}{\bar{z}}
\newcommand{\zbyneg}{\ensuremath{\zbar_{(+-)}}}
\newcommand{\zbxneg}{\ensuremath{\zbar_{(-+)}}}
\newcommand{\zbell}{\ensuremath{\zbar_{(e)}}}
\newcommand{\eps}{\varepsilon}
\renewcommand{\t}[1]{\tilde{#1}}
\newcommand{\doubleslash}{/\hspace{-3pt}/}
\newcommand{\Ep}{\textit{E}-piece}
\DeclareMathOperator{\tr}{tr}

\DeclareMathOperator{\Out}{Out}
\newlength{\BS}
\newcommand{\nl}{\\*[\BS]}
\theoremstyle{plain}
\newtheorem{thm}{Theorem}
\newtheorem*{mainA}{Theorem A}
\newtheorem*{mainB}{Theorem B}
\newtheorem{lemma}[thm]{Lemma}
\newtheorem{prop}[thm]{Proposition}
\newtheorem{cor}[thm]{Corollary}
\newtheorem*{conjecture}{Conjecture}
\newtheorem*{defin}{Definition}

\newtheorem*{claim}{Claim}
\theoremstyle{definition}
\newtheorem*{remark}{Remark}
\numberwithin{equation}{section}
\numberwithin{thm}{section}
\newcommand{\Hom}{\ensuremath{\mop{Hom}}}
\newcommand{\Aut}{\ensuremath{\mop{Aut}}}
\newcommand{\lrarw}{\longrightarrow}
\newcommand{\Tei}{{\mathfrak T}}
\newcommand{\Ht}{\mathbf{H}^2}
\newcommand{\slmr}{\ensuremath{\mop{SL}_{-}(2, \mathbb{R})}}
\newcommand{\slpmr}{\ensuremath{\mop{SL}_{\pm}(2, \mathbb{R})}}
\newcommand{\islr}{\ensuremath{\mop{ISL}(2, \mathbb{R})}}
\newcommand{\blm}{\ensuremath{B_{\Lambda}}}
\newcommand{\blmd}{\ensuremath{\blm(u)}{}}
\newcommand{\hcz}{\ensuremath{h_{c,z}}}
\newcommand{\hczz}{\ensuremath{h_{c,z_0}}}
\newcommand{\lxz}{\ensuremath{L_{x-x_0}}}
\newcommand{\lyz}{\ensuremath{L_{y-y_0}}}
\newcommand{\lzz}{\ensuremath{L_{z-z_0}}}
\newcommand{\scz}{\ensuremath{\Sigma_0^c}}
\newcommand{\ecb}{\ensuremath{\overline{\E{E}}_c}}
\newcommand{\mcv}{\ensuremath{\mathfrak{m}_c(v)}}
\newcommand{\tpp}{\texttt{(+,+)}}
\newcommand{\tpm}{\texttt{(+,-)}}
\newcommand{\tmp}{\texttt{(-,+)}}
\newcommand{\tzz}{\texttt{(0$|$0)}}
\newcommand{\rpt}{{\R\mathbb{P}}^2}
\newcommand{\linf}{l_{\infty}}
\newcommand{\pinf}{P^{\infty}}
\floatstyle{boxed}
\restylefloat{figure}

\begin{document}
\title[Action of the Modular Group]
{Dynamics of  the Automorphism Group of the
$GL(2,\mathbb{R})$-Characters of a Once-puncutred Torus}
\author{William Goldman}
\author{George Stantchev}
\address{Department of Mathematics \\
University of Maryland \\
College Park, MD 20742}
\date{May 23, 2003}
\begin{abstract}
Let $\pi$ be a free group of rank $2$. Its outer automorphism group
$\Out(\pi)$ acts on the space of equivalence classes of
representations $\rho \in \hom{\pi}{\slc}$. Let $\slmr$ denote the
subset of \glr{} consisting of matrices of determinant $-1$ and let
$\islr$ denote the subgroup $\slr \amalg i \slmr \subset \slc$.
The representation space $\hom{\pi}{\islr}$ has four 
connected components, 
three of which consist of representations 
that send at least one generator of $\pi$ to $i\slmr$. 
We investigate the dynamics of the $\Out(\pi)$-action on these  
components. 
\par
The group $\Out(\pi)$ is commensurable with the group $\Gamma$ of
automorphisms of the polynomial
\[
\kappa(x,y,z) = -x^2 - y^2 + z^2 + xyz -2 
\]
We show that for $-14 < c < 2$, the action of  $\Gamma$ is ergodic on
$\kappa\inv(c)$. For $c< -14$, the group $\Gamma$ acts 
properly and freely on an open subset 
$\Omega^M_c \subset \kappa\inv(c)$ and acts ergodically on
the complement of $\Omega^M_c$. We construct an algorithm 
which determines, in polynomial time,  if a point $(x,y,z)\in\R^3$ is 
$\Gamma$-equivalent to a point in $\Omega^M_c$ or in its complement.
\par
Conjugacy classes of $\islr$-representations identify with $\R^3$ via
the appropriate restriction of the character map
\begin{align*}
\chi:\,\hom{\pi}{\slc} & \longrightarrow \C^3 \\
\rho  & \longmapsto 
\bmatrix \xi(\rho) \\ \eta(\rho)  \\ \zeta(\rho) \endbmatrix  = 
\bmatrix \tr (\rho (X))\\
\tr (\rho (Y))\\ \tr (\rho (XY))\endbmatrix
\end{align*}
where $X$ and $Y$ are the generators of $\pi$. 
Corresponding to the Fricke spaces of the once-punctures Klein bottle
and the once-punctured M\"obius band are $\Gamma$-invariant open
subsets $\Omega^K$ and $\Omega^M$ respectively. 
We give an explicit parametrization of $\Omega^K$ and $\Omega^M$ as
subsets of $\R^3$ and we show that $\Omega^M\cap \kappa\inv(c)\neq
\varnothing$ if and only if $c<-14$, while $\Omega^K\cap \kappa\inv(c)\neq
\varnothing$ if and only if $c>6$. 

\end{abstract}

\maketitle
\tableofcontents
\listoffigures


\section{Introduction}
Let $\pi$ be a free group of rank $2$. Its outer automorphism group
$\Out(\pi)$ acts on the space of equivalence classes of
representations $\rho \in \hom{\pi}{\slc}$. Let $\slmr$ denote the
subset of \glr{} consisting of matrices of determinant $-1$, and let
\[
\slpmr = \{A \in \glr \mid \mop{det}(A) = \pm 1\}
\]
The group $\slpmr$ is isomorphic to 
\[
\islr = \slr \amalg i \slmr 
\]
and in this context we identify the two as subgroups of $\slc$.
The representation space $\E{R}=\hom{\pi}{\islr}$ has four connected
components indexed by the 
elements of \linebreak[1] 
\mbox{$\mop{H}^1(\pi, \Z_2) \cong \Z_2\times\Z_2$}.
The three non-zero elements of $\Z_2\times\Z_2$ correspond to the
components of 
$\E{R}$ consisting of representations that send at least one generator
of $\pi$ to $i\slmr$. We investigate the dynamics of the
$\Out(\pi)$-action on these 
components. The action of $\Out(\pi)$ on the component of
$\slr$-representations has been recently studied by
Goldman~\cite{Gold:pTorus}
\par By a theorem of Fricke~\cite{FrickeKlein}, the moduli space of
$\slc$-representations naturally identifies with affine $3$-space
$\C^3$ via the character map
\begin{align*}
\chi:\,\hom{\pi}{\slc} & \longrightarrow \C^3 \\
\rho  & \longmapsto 
\bmatrix \xi(\rho) \\ \eta(\rho)  \\ \zeta(\rho) \endbmatrix  = 
\bmatrix \tr (\rho (X))\\
\tr (\rho (Y))\\ \tr (\rho (XY))\endbmatrix
\end{align*}
where $X$ and $Y$ are the generators of $\pi$. Let $[X,Y]$ be the
commutator of $X$ and $Y$. In
terms of the coordinate functions $\xi$, $\eta$ and $\zeta$,  the
trace $\tr([X,Y])$ is given by the polynomial 
\[
\kappa(\xi, \eta, \zeta) := \xi^2 + \eta ^2 + \zeta^2 - \xi\eta\zeta -2
\]
which is preserved under the action of $\Out(\pi)$. Moreover, the
action of $\Out(\pi)$ on $\C^3$ is commensurable with the action
of the group $\Gamma$ of polynomial automorphisms of $\C^3$ which
preserve $\kappa$ (Horowitz~\cite{Horowitz}). Note that $\Gamma$ is
a finite extension of the modular group and is isomorphic to
\[
\pglz \ltimes (\Z/2\oplus\Z/2)
\]
\par 
Let $\rep{1}{1}$ be the component of $\Hom(\pi, \islr)$ consisting 
of representations that send both $X$ and $Y$ to $i\slmr$. 
The restriction $\chi_{11}$ of the character map to $\rep{1}{1}$ is a
 surjection onto $\cv{1}{1} = i\R\times i\R \times \R$. The latter
is isomorphic to $\R^3$ and 
the restriction of $\kappa$ to $\cv{1}{1}$ induces a polynomial
\[
\kappa_{11}(x,y,z) := \kappa(ix, iy, z) = -x^2 - y^2 + z^2 + xyz -2 
\]
where $x$, $y$ and $z$ are the standard coordinate functions on
$\R^3$. If $u\in \cv{1}{1}$ is such that $k(u) \neq 2$, the fiber 
$\chi_{11}\inv(u)$ is an $\slpmr$-conjugacy class of irreducible
representations in $\rep{1}{1}$. In this context, $\cv{1}{1}$
identifies with a component of the $\slpmr$-character variety of
$\pi$. 
\begin{mainA}\label{T:ergintro}
Let $\kappa_{11}(x,y,z) = -x^2 -y^2 +z^2 + xyz -2$ and let $c\in
\R$. Let $\Gamma$ be the automorphism group of $k_{11}$. Then 
\begin{itemize}
\item For $-14\leq c < 2$, the group $\Gamma$ acts ergodically on
$\kappa_{11}\inv(c)$. 
\item For $c< -14$, the group $\Gamma$ acts properly and freely on an
open subset $\Omega^M_c \subset \kappa_{11}\inv(c)$ and acts ergodically on
the complement of $\Omega^M_c$. 
\end{itemize}
\end{mainA}
Since $\pi$ is a free group, representations in $\hom{\pi}{\slpmr}$,
or equivalently in \linebreak[4] $\hom{\pi}{\islr}$,
can be realized as lifts of representations in $\hom{\pi}{\pglr}$ and
thus can be interpreted geometrically
via the identification of \linebreak[1] $\pglr$ with the full 
isometry group of hyperbolic $2$-space $\Ht$. More precisely, let $S$ be 
a surface  with fundamental group $\pi_1(S)$. Let $G$ be a semisimple Lie
group. Then $\hom{\pi_1(S)}{G}$  is an analytic variety upon which G
acts by conjugation. Let $\hom{\pi_1(S)}{G}/G$ be the orbit space. 
The $G$-orbits parametrize equivalence classes of flat principal
$G$-bundles over $S$. If $X$ is a space upon which $G$ acts,
$\hom{\pi_1(S)}{G}/G$ is the 
deformation space of flat $(G,X)$-bundles over $S$. In this context,
\mbox{$\hom{\pi}{\pglr}/{\pglr}$} identifies with the deformation space of flat
$\Ht$-bundles over a surface $S$ whose fundamental group is free of
rank $2$. 
\par When $\rho\in \hom{\pi}{\pglr}$ is a discrete embedding, the
holonomy group $\rho(\pi)$ acts properly discontinuously on the fiber
$\Ht$. The quotient $\Ht/\rho(\pi)$ is homotopy-equivalent to $S$
and affords a hyperbolic structure induced by that of $\Ht$. If the
quotient is also diffeomorphic to $S$ we call $\rho$ a \emph{discrete
$S$-embedding}. The set  $\Omega^S$ of conjugacy classes of discrete
$S$-embeddings is open in $\hom{\pi}{\pglr}/{\pglr}$ and parametrizes complete 
hyperbolic structures on $S$ marked with respect to a fixed set of
generators of $\pi$. We call $\Omega^S$ the \emph{Fricke space} of $S$. 
In a certain sense $\Omega^S$ is  a generalization of the Teichm\"uller
space of an orientable closed surface. 
\par Discrete embeddings inside  $\rep{1}{1}$ give rise to
non-orientable surfaces and since $\pi$ is free of rank $2$ the only
possibilities are the once-punctured%
\footnote{By ``punctured'' in this context we mean a surface with
the interiors of one or more disjoint disks removed, i.e.
a surface with boundary}
M\"obius band $M$ (equivalently, the twice-punctured projective
plane), and the once-punctured Klein bottle $K$. 
Their respective Fricke spaces $\Omega^M$ and $\Omega^K$ can
be parametrized as subsets of $\R^3$. 
\begin{mainB}\label{T:frickeintro}
Let $\Gamma$ be the group of automorphisms of the polynomial
\[
\kappa_{11}(x,y,z) = -x^2 - y^2 + z^2 - xyz -2
\]
\begin{enumerate}
\item
Let $\oM$ be the region in $\R^3$ defined by the inequalities
\begin{align*} 
xy + z &> 2 \\
z & < - 2 
\end{align*}
Then the Fricke space  of the once-punctured M\"obius band
$M$ identifies with
\begin{equation*}
\Omega^M = \coprod_{\gamma\in\Gamma/\Gamma_M}\gamma \oM
\end{equation*}
where $\Gamma_M\subset\Gamma$ is the stabilizer of $\oM$, 
and $\Gamma/\Gamma_M$ denotes the  coset space of $\Gamma_M$. 
\item Let $\Omega_0^K$ be the region in $\R^3$ defined by the
inequality
\begin{equation*}
x^2 + y^2 - xyz +4 < 0
\end{equation*}
Then the Fricke space of the once-punctured Klein bottle
$K$ identifies with
\begin{equation*}
\Omega^K = \coprod_{\gamma\in\Gamma/\Gamma_K}\gamma \Omega_0^K
\end{equation*}
where $\Gamma_K\subset\Gamma$ is the stabilizer of $\Omega_0^K$, and
$\Gamma/\Gamma_K$ denotes the  coset space of $\Gamma_K$. 
\end{enumerate}
\begin{remark}
The subgroups $\Gamma_M$ and $\Gamma_K$ correspond to the mapping
class groups of $M$ and $K$ 
respectively. 
\end{remark}
\end{mainB}
The level sets  $\kappa_{11}\inv(c)$ intersect $\Omega^M$ if and only if
$c<-14$,  and they intersect $\Omega^K$ if and only if $c>6$. Assume $c<-14$,
and let $\Omega_c^M = \kappa_{11}\inv(c)\cap \Omega^M$. This is precisely the
region mentioned in the second part of Theorem~A. 
In other words, on the Fricke space of the once-punctured M\"obius band $M$
the action of $\Gamma$ is wandering; on the outside, the action is ergodic along the
level sets  $\kappa_{11}\inv(c)$ for each $c<2$.  Similar statements
hold for the Fricke spaces of $M$ sitting inside the other
$\slpmr$-moduli-space components that correspond to non-zero classes in
$\Z_2\times\Z_2$.
\begin{conjecture} Let $\Omega^K$ be the Fricke space of the once-punctured
Klein bottle. Let $\Omega_c^K = \kappa_{11}\inv(c)\cap \Omega^K$
\begin{itemize}
\item The action of $\Gamma$ on $\Omega^K$ is wandering. 
\item For each $c>2$, the action of $\Gamma$ on the set
$\kappa_{11}\inv(c) - \Omega_c^K$ is ergodic. 
\end{itemize}
\end{conjecture}

We are grateful to Joan Birman, Nikolai Ivanov, Misha Kapovich, John
Millson, Walter Neumann and Scott Wolpert for helpful discussions on
this material.

\section{Background and Motivation}
In this part we provide some relevant background material and 
interpret it in our context. We also introduce new notation
and terminology that will be used in the subsequent 
exposition.
\subsection{Algebraic Generalities}
If $A$,$B$ are groups, then $\Hom(A,B)$ denotes the set of homomorphisms
(representations) $A \lrarw B$. Let \Aut(A) denote the group of all
automorphisms of $A$. There is an action of the group $\Aut(A) \times
\Aut(B)$ on $\Hom(A,B)$ defined by
\begin{equation} \label{E:actionAut}
\begin{aligned}
        (\Aut(A) \times \Aut(B)) \times \Hom(A,B) & \lrarw \Hom(A,B)\\
 ((\alpha, \beta), \phi) & \longmapsto \beta \circ \phi \circ
 \alpha^{-1}
\end{aligned}
\end{equation}
In particular, the action of $B$ on itself by conjugation embeds $B$
in its automorphism group and thus induces an action:
\[
\begin{aligned}
        B \times \Hom(A,B) & \lrarw \Hom(A,B)\\
 ( b , \phi) & \longmapsto \iota_b \circ \phi 
\end{aligned}
\]
where $\iota_b: h\mapsto bhb^{-1}$ is the inner automorphism of
$B$ defined by conjugation by 
$b \in B$. The orbit space, $\Hom(A,B)/B$ is the set of conjugacy
classes of representations in $\Hom(A,B)$ and the action
(\ref{E:actionAut}) descends to an action of $\Aut(A)$ on $\Hom(A,B)/B$.
Let $\mop{Inn}(A)$ denote the (normal) subgroup of $\Aut(A)$ consisting
of \emph{inner automorphisms}.
Since $\mop{Inn}(A)$ preserves the conjugacy class of a
representation, it acts trivially on $\Hom(A,B)/B$ and thus the action
of $\Aut(A)$ factors through the action of the \emph{outer
automorphism group} 
\[
\Out(A):=\Aut(A)/\mop{Inn}(A)
\]
\subsection{Geometric Motivation}
When $A$ is a discrete group and $B$ a Lie group, the representation
space $\Hom(A,B)$ can have special geometric significance. In
particular, it can be interpreted as the moduli space of 
flat bundles over manifolds. 
Consider, for example, a manifold $M$ with fundamental group
$\pi=\pi_1(M,x_0)$ and a Lie group $G$ which acts on a space $X$.
Then the space of \emph{flat $(G,X)$-bundles over $M$} 
can be identified with
\hpgmg. In the case when $M$ is a closed surface, and $G$ is the group
of orientation preserving isometries of the hyperbolic plane, the
subset of $\hpgmg$ consisting of equivalence classes of discrete
embeddings of $\pi$, identifies with the Teichm\"uller space $\Tei$ of
$M$. Moreover, $\Out(\pi)$ is isomorphic to the mapping class group
$\pi_0\mop{Diff}(M)$ of $M$ whose action on $\Tei$ is well known to be
\emph{properly discontinuous}.
\subsection{The Structure of $\hpg$}
Whenever $\pi$ is a finitely generated group, the space \hpg{} inherits
a natural (Hausdorff) topology from that of $G$. Namely, if 
$\{\gamma_1,\dots, \gamma_n\}$ is a set of generators of $\pi$ then
the evaluation map
\[
\begin{aligned}
\hpg & \lrarw G^n \\
\rho & \longmapsto (\rho(\gamma_1), \dots, \rho(\gamma_n))
\end{aligned}
\]
is an embedding which induces a topology on \hpg{} that is independent
on the choice of generators. Furthermore, if $G$ is an algebraic
group, then there is an induced \emph{algebraic structure} on \hpg{}, and
we refer to the resulting object as the \emph{variety of
representations}. 
\nl
However, the topology of the space $\hpgmg$ inherited from \hpg{}
could be rather pathological (e.g. may not be even Hausdorff). 
When $G$ is a reductive linear algebraic group, i.e a subgroup of 
$\o{GL}(n,\mathbb{R})$, however, one can consider the
algebraic-geometric quotient $\hpg\doubleslash G$, which is defined as the variety
whose coordinate ring is precisely the ring of $G$-invariant functions
on $\hpg{}$. Equivalently, it is the space of \emph{semi-stable orbits}.
In this case, the orbit of a representation $\rho$ is semi-stable if 
$\rho$ is a completely reducible representation. We call $\hpg\doubleslash G$ the
\emph{character variety} of $\pi$. 
\subsection[Obstruction Classes and Components of $\hpg$]%
{Obstruction Classes and Components \\ of $\hpg$ }
\label{SS:Stiefel}
Let $S$ be a connected surface and let $\pi$ denote its fundamental
group. Let $G$ be a real algebraic Lie group. An element $\rho$
in the space $\hpg{}$  determines a flat prinipal $G$-bundle over
$S$. This bundle gives rise to the obstruction classes
in $H^q(S, \pi_{q-1}(G))$ and thus induces an
\emph{obstruction class map}
\[
o_q: \hpg \longrightarrow H^q(S, \pi_{q-1}(G))
\]
In particular, let $G_0$ be the identity component for $G$ and
$G/G_0= \pi_0(G)$ be the group of components of $G$. Then
$o_1(\rho) \in H^1(S, \pi_0(G))$, which via the Hurewicz
isomorphism is $\hom{\pi}{G/G_0}$. Thus, $o_1(\rho)$ is just the
composite 
\[
\pi \xrightarrow[]{\rho} G \longrightarrow G/G_0.
\]
The map $o_1$ is continuous, hence, constant on each connected
component. 
For example if $G=\pglr \simeq \mop{SO}(2,1)$, then
$\pi_0(G) = \mbb{Z}_2$ and 
$o_1(\rho)\in H^1(S, \mbb{Z}_2)$. In that case the action of $\pglr$
on hyperbolic $2$-space $\Ht$ gives rise to an associated 
circle bundle (with fiber $\partial\Ht = \mathbb{RP}^1$), whose first 
Stiefel-Whitney class corresponds to $o_1(\rho)$ (see
Steenrod~\cite{Steen}, \S$38$; Goldman~\cite{TopComp},
\S$2$). 
If $n>0$ then $\pi$ is a free group of rank
$r=1-\chi(S)$, where $\chi(S)$ is the Euler characteristic of
$S$. Thus $H^1(S, \mbb{Z}_2) \simeq \mbb{Z}_2^r$, which is a free
$\mbb{Z}_2$-module. In this case $\hom{\pi}{G}\cong G^r$ and the
fibers of $o_1$ identify with the connected components of
$\hom{\pi}{G}$. In this context we shall refer to  $o_1$ as  the
\emph{first Stiefel-Whitney class map}.  
\subsection{$\slc$-character Varieties}
The example in the last paragraph is of particular geometric
interest. Consider a finitely generated group $\pi$ and the
representation space $\hom{\pi}{\slc}$. Since \slc{} acts as a 
group of symmetries of hyperbolic $3$-space $\mathbf{H}^3$, 
the orbit space $\hom{\pi}{\slc}/\slc$ contains a subset that
parametrizes 
equivalence classes  of hyperbolic structures on $3$-manifolds with 
fundamental group isomorphic to $\pi$. 
\par

The character variety
\[
\E{X} = \hom{\pi}{\slc}\doubleslash \slc
\]
admits an embedding of $\E{X}$ as an algebraic subset of affine space. 
In particular, 
traces of a finite generating set of $\pi$. 
define coordinates on $\E(X)$ (Procesi~\cite{Procesi}).
We consider the case when $\pi$ is free of rank $2$:

\begin{thm}[Fricke-Klein]
Let $X$, and $Y$ be the generators of $\pi$ and let $G = \slc$. Then
the character map: 
\begin{align*}
\chi: \hom{\pi}{G}\doubleslash G  & \longrightarrow \C^3 \\
[\rho]  & \longmapsto 
\bmatrix \xi(\rho) \\ \eta(\rho)  \\ \zeta(\rho) \endbmatrix  = 
\bmatrix \tr (\rho (X))\\
\tr (\rho (Y))\\ \tr (\rho (XY))\endbmatrix
\end{align*}
is an isomorphism. (Compare the discussion in Goldman~\cite{TopComp},
4.1, \cite{Erg},\S\S4--5 and ~\cite{Expo}.) Thus the traces of $X, Y,
XY$ parametrize $\hpg\doubleslash G$ as the affine space $\C^3$. In particular,
if $w(X,Y)$ is any word in $X$  and $Y$, then $\mop{tr}(\rho(w(X,Y)))$
is expressed as a polynomial $f_w$ in $\xi,\eta, \zeta$.
\end{thm}
\par

Let \E{H} be a copy of the hyperbolic plane $\Ht$ sitting inside
$\mathbf{H}^3$. The stabilizer of \E{H} in $\slc$ is isomorphic to
\pglr{} and in this way \pglr{} identifies with the full group of
isometries of $\Ht$. Naturally, we are interested in the orbits of
the action of $\pglr$ by conjugation on $\hom{\pi}{\pglr}$. 
However, in order to use the trace parametrization of the orbit space 
provided by the theorem of Fricke-Klein, we need to work with 
the relevant representation space, which in this case is
$\hom{\pi}{\slpmr}$.  The subgroup 
\[
\slpmr := \{A \in \glr\mid \mop{det}(A) = \pm 1\}
\]
of $\glr$ is a double cover of $\pglr$.
Since $\pi$ is a free group, any representation $\rho\in
\hom{\pi}{\pglr}$ lifts to a representation in $\hom{\pi}{\slpmr}$. 
Thus a conjugacy class of representations in $\hom{\pi}{\pglr}$
corresponds to a conjugacy class of representations in
$\hom{\pi}{\slpmr}$ together with a choice of a lift. 
\subsection{$\slpmr$-character Varieties}
Let $\slmr$ denote the
subset of \glr{} consisting of matrices of determinant $-1$.
The group $\slpmr$ is isomorphic to 
\[
\islr = \slr \amalg i \slmr 
\]
and in this context we identify the two as subgroups of $G=\slc$.
In view of the discussion in section~\ref{SS:Stiefel},
the representation space 
\[
\E{R}=\hom{\pi}{\slpmr}\cong\hom{\pi}{\islr} \cong \islr \times \islr
\] 
has four connected components indexed by the 
elements of $\Z_2\times\Z_2$. 
Let $G_0 = \slr$ and $G_1 = i\slmr$.
Then the correspondence defined by the Stiefel-Whitney class map is
\[
\rep{j}{k} = G_j\times G_k \longmapsto (j,k)\in \Z_2\times\Z_2
\]
where $j,k \in \{0,1\}$. The following Proposition shows that the
the restriction of the character map $\chi$ to  $\E{R}$ can be
used to parametrize conjugacy classes of $\islr$ representations
(compare Xia~\cite{Xia}, pp.~10-13). 
\begin{prop}\label{P:xiaprop}
Let $\chi_{jk}$  be the restriction of $\chi$  to $\rep{j}{k}$.
Assume $(j,k) \neq (0,0)$. Then 
\begin{enumerate}
\item $\chi_{jk}$ is surjective
\item The image of $\chi_{jk}$ is
\[
\cv{j}{k} := \chi(\rep{j}{k}) = \left\{
\begin{array}{lll}
\R\times i\R \times i\R  & \text{if} & (j,k) = (0,1) \\
i\R\times \R \times i\R  & \text{if} & (j,k) = (1,0) \\
i\R\times i\R \times \R  & \text{if} & (j,k) = (1,1) \\
\end{array}
\right.
\]
\item Let $\rho\in \E{R}$. Then 
\begin{equation}\label{E:kappa}
\kappa\circ\chi(\rho) = \left\{
\begin{array}{lll}
\phantom{-}x^2 -y^2 -z^2 + xyz -2 & \text{if} & \rho \in \rep{0}{1} \\
-x^2 +y^2 -z^2 + xyz -2 & \text{if} & \rho \in \rep{1}{0} \\
-x^2 -y^2 +z^2 + xyz -2 & \text{if} & \rho \in \rep{1}{1} \\
\end{array}
\right.
\end{equation}
\item Let $u\in \cv{j}{k}$ be such that $k(u) \neq 2$. Then $\islr$
acts transitively on $\chi_{jk}\inv(u)$. 
\end{enumerate}
\begin{proof}
We refer the proof to Xia~\cite{Xia}, Proposition $12$ and $13$. 
\end{proof}
\end{prop}
By a result of Culler and Shalen, a representation in $\rho
\in\E{R}$  is reducible if and only if $\kappa(\chi(\rho))=2$. 
Thus  Proposition~\ref{P:xiaprop} implies that if $u\in \cv{j}{k}$ is not in 
$\kappa\inv(2)$, then the fiber  
$\chi_{jk}\inv(u)$ is an $\islr$-conjugacy class of irreducible
representations in $\rep{j}{k}$. In that sense $\cv{j}{k}$
identifies with a component of the $\islr$-character variety of
$\pi$. 
\par
Let $\kappa_{jk}$ be the restriction of $\kappa$ to $\cv{j}{k}$. 
Since $\cv{j}{k}$ is naturally isomorphic to $\R^3$ we will 
in fact consider the fibers of $\kappa_{jk}$ as subsets of $\R^3$. 
Furthermore, whenever it is clear from the context, we will drop the
subscript and use $\kappa$ to denote the relevant  polynomial as
prescribed in~(\ref{E:kappa}).
\subsection{Surfaces Whose Fundamental Group is Free of Rank $2$}
\label{SS:surf-rank-two}
\begin{figure}[htbp]
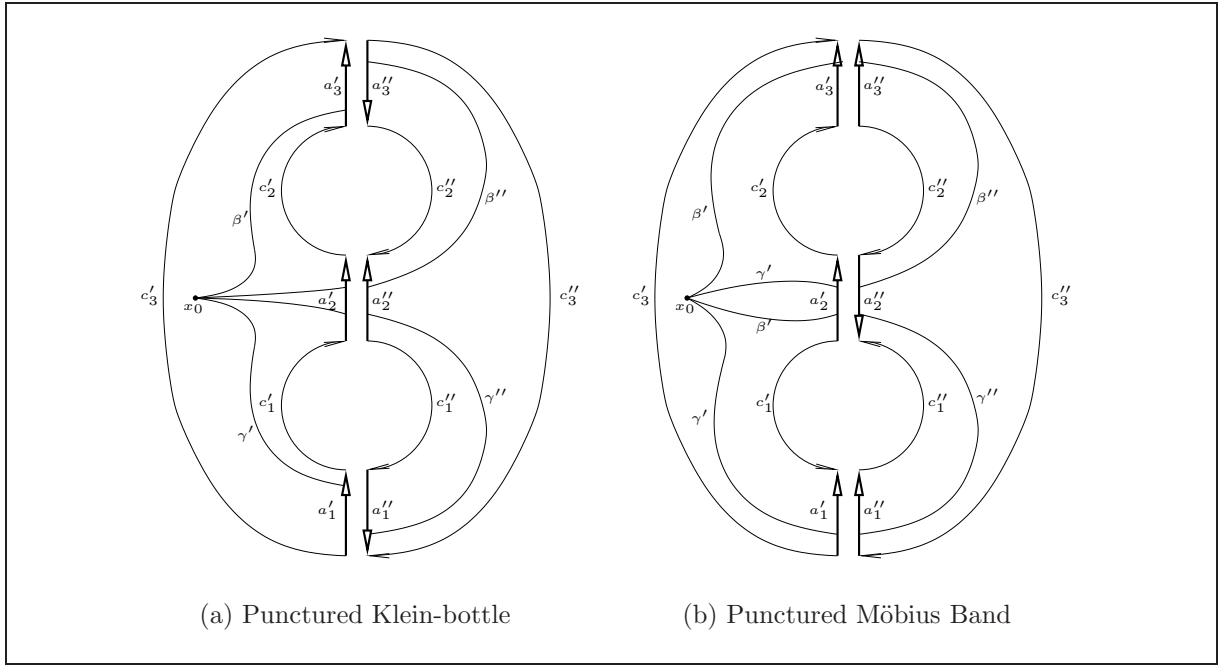

\centering
\subfigure[Punctured Klein-bottle]{\label{F:pk-bottle}
\input{pK-bottle.pstex_t}
}
\quad
\subfigure[Punctured M\"obius Band]{\label{F:pM-band}
\input{pM-band.pstex_t}
}
\caption{Constructing non-orientable surfaces from E-pieces}\label{F:pmb-pkb}
\end{figure}

Let $S$ be a surface whose fundamental group $\pi$ is free of rank
$2$. Then $S$ is precisely one of the following: 
a punctured\footnote{In this context a ``punctured surface'' shall
mean a surface with a topological disk removed} 
torus $S_{1,1}$, a pair of pants $S_{0,3}$, a punctured Klein bottle
$S^{\sharp}_{1,1}$, or a 
punctured M\"obius band $S^{\sharp}_{0,2}$.
These are all connected surfaces with Euler characteristic $-1$. 
\begin{remark} The notation used for referencing these surfaces
is based on the genus and the number of boundary components of each
surface. Thus $S_{g,n}$ (respectively $S_{g,n}^{\sharp}$) denotes an
orientable (respectively non-orientable) surface of genus $g$ and 
$n$ boundary components. For non-orientable surfaces, however, the
notion of genus is not very consistently used throughout the literature.
Some authors define genus as the number of copies of the projective
plane in the cross-cap decomposition of the surface. 
In that context, the genus of the Klein bottle would be $2$, and
the genus of the M\"obius band would be $1$. 
\par We adhere to the definition that is consistent with the formula
for the Euler characteristic of a non-orientable surface $S$
\[
\chi(S) = 1 - \beta_1(S) - n
\]
where $\beta_1(S)$ is the first Betti number and $n$ is the number of
boundary components. We take $\beta_1(S)$ as a definition for  the
genus of a non-orientable surface $S$ (see Moise~\cite{Moise}, 
\S~$22$).  
\end{remark}
Next, we  show how the geometry of each  surface is related to the 
presentation of $\pi$. 
\begin{defin}
An \Ep{} is a right-angled hyperbolic hexagon whose boundary
is a piecewise-geodesic closed curve without self-intersections. Each
edge of an \Ep{} is given the orientation of the underlying
geodesic segment. 
\end{defin}
Each surface of Euler characteristic $-1$ can be realized by pasting
two \Ep{}s subject to a certain gluing scheme. Figure~\ref{F:pmb-pkb}
illustrates this idea for the non-orientable cases.
Let $\rho \in \hom{\pi}{\pglr}$ be a discrete $S$-embedding. The union of two
\Ep{}s $E'$ and $E''$ is a fundamental  
region for $\rho(\pi)$ acting on $\Ht$ with quotient $S$. From the
gluing diagrams we obtain a geometric presentation of $\pi$ in
each case of a surface of Euler characteristic $-1$. For example,
suppose $S=\PKB$ is the punctured Klein-bottle. Then, 
$S = E' + E''$ modulo the gluing pattern shown on
Fig.~\ref{F:pk-bottle}. In the quotient, $\beta' * \beta'' = \beta$,
and $\gamma' * \gamma'' = \gamma$  become the two generators of $\pi$,
and the boundary geodesic that represents the puncture is $\delta = 
c_2'*c_3''*c_1'*c_1''*c_3'*c_2''$. Here the ``$*$'' operator indicates
the obvious left-to-right curve concatenation with respect to a suitable
parametrization on each geodesic segment. Thus, for a fixed base point
$x_0$, we obtain a presentation of $\pi \equiv \pi_1(S, x_0)$ as follows:
\begin{equation}\label{E:pKb-11pres}
	\pi = \langle \beta, \gamma, \delta \mid \gamma \beta^2 \gamma =
	\delta \rangle
\end{equation}
Similarly, from the diagram in Fig.~\ref{F:pM-band} we obtain the following
presentation of the fundamental group $\pi \equiv
\pi_1(\PMB, x_0)$ of the punctured M\"obius band:
\begin{equation}
	\pi = \langle \beta, \gamma, \delta_1, \delta_2 \mid \delta_1 =
	\beta \gamma, \delta_2 = \beta\gamma\inv \rangle
\end{equation}
where again $\beta = \beta'*\beta''$ and $\gamma = \gamma'*\gamma''$
are the two generators of $\pi$, while
$\delta_1 = c_3'*c_3''$ and $\delta_2 = c_2'*c_2''*c_1'*c_1''$ 
represent the two boundary components of $\PMB$. 
\subsection{Mapping Class Group, and the Structure of $\Out(\pi)$} 
For a compact connected surface $S$, the Mapping Class Group of $S$ is
defined as the group $\pi_0(\mop{Homeo}(S,\partial S))$ of isotopy
classes of homeomorphisms on $S$. There is a well-defined
homomorphism:
\begin{equation}\label{Eq:Nielsen}
N:\MCG \longrightarrow \Out(\pi) \cong \Aut(\pi)/\mop{Inn}(\pi).
\end{equation}
If $S$ is a closed surface then by Dehn (unpublished) and Nielsen~
\cite{Nielsen1927}, $N$ is an isomorphism. 
When $\partial S\neq\emptyset$, then each component $\partial_i S$
determines a conjugacy class $C_i$ of elements of $\pi_1(S)$
and the image of $N$ consists of elements of $\Out(\pi)$ represented by
automorphisms which preserve each $C_i$. Another theorem of
Nielsen~\cite{Nielsen}, implies that when $S$ is a punctured torus,
$N$ is also an isomorphism.
\nl
The action of $\Out(\pi)$ on homology $H_1(S;\Z)\cong \Z^2$ defines a homomorphism
\begin{equation}\label{eq:homology}
h:\Out(\pi) \longrightarrow \glz.
\end{equation}
which, again by Nielsen~\cite{Nielsen}, is an isomorphism (see also
Magnus-Karrass-Solitar~\cite{Magnus-KS}, \S~3.5, Corollary N4).
Thus for a free group $\pi$ of rank $2$, $\Out(\pi)$ is isomorphic to $\glz$.

\subsection{The modular group}
In general the group $\Aut(\pi)$ acts on the $\slc$ character variety
$\E{X}$ by:  
\begin{equation*}
\phi_*([\rho]) = [\rho\circ\phi^{-1}]
\end{equation*}
Let $\gamma\in \pi$ and let $\iota_\gamma$ denote conjugation by
$\gamma$. Since
\begin{equation*}
\rho\circ\iota_\gamma = \iota_{\rho(\gamma)}\circ\rho
\end{equation*}
the subgroup $\mop{Inn}(\pi)$ acts trivially. Thus $\Out(\pi)$ acts on
$\C^3$
and since an automorphism $\phi$ of $\pi$ is determined by 
\begin{equation*}
(\phi(X),\phi(Y))= (w_1(X,Y),w_2(X,Y)),  
\end{equation*}
for some words $w_1$, $w_2$ in the generators, 
the action of $\phi$ on
$\C^3$ is given by a triple of polynomials
\begin{equation*}
(f_{w_1}(\xi,\eta,\zeta),f_{w_2}(\xi,\eta,\zeta),
f_{w_1w_2}(\xi,\eta,\zeta))
\end{equation*}
Hence $\Out(\pi)$ acts on $\C^3$ by {\em polynomial automorphisms.\/}
Nielsen's theorem (see Magnus-Karras-Solitar~\cite{Magnus-KS}, Theorem 3.9)
implies that any such automorphism preserves
$\kappa:\C^3\longrightarrow\C$, 
that is 
\begin{equation*}
\kappa\left((f_{w_1}(\xi,\eta,\zeta),f_{w_2}(\xi,\eta,\zeta),
f_{w_1w_2}(\xi,\eta,\zeta))\right) 
= \kappa(\xi,\eta,\zeta).
\end{equation*}
Horowitz\cite{Horowitz} determined the group $\Aut(\C^3,\kappa)$ of
polynomial mappings  $\C^3\longrightarrow\C^3$ preserving $\kappa$.
Let $\mathfrak{S}_3$ be the symmetric group consisting of permutations
of the coordinates $\xi, \eta, \zeta$. Horowitz
proved that the automorphism group of 
$(\C^3,\kappa)$ is generated by the linear automorphism group
\begin{equation*}
\Aut(\C^3, \kappa)\cap \mop{GL}(3,\C) = 
(\Z_2\oplus\Z_2)\rtimes\mathfrak{S}_3
\end{equation*}
and the {\em quadratic reflection:\/}
\begin{equation*}
\bmatrix \xi \\ \eta \\ \zeta \endbmatrix
\longrightarrow \bmatrix \xi \\ \eta \\ \xi\eta - \zeta \endbmatrix,
\end{equation*}
This group is commensurable with $\Out(\pi)$. The factor
$\Z_2\oplus\Z_2$ corresponds to \emph{sign change automorphisms} of
$\E{R}$. In particular, the three non-trivial elements $(0,1),
(1,0),(1,1)$ act on representations by $\sigma_{10}$, $\sigma_{01}$,
$\sigma_{11}$ respectively:
\begin{equation*}
\sigma_{01}\cdot\rho: 
\begin{cases}
X & \longmapsto \rho(X) \\
Y & \longmapsto -\rho(Y) \end{cases}
\end{equation*}
\begin{equation*}
\sigma_{10}\cdot\rho: 
\begin{cases}
X & \longmapsto -\rho(X) \\
Y & \longmapsto \rho(Y) \end{cases}
\end{equation*}
\begin{equation*}
\sigma_{11}\cdot\rho: 
\begin{cases}
X & \longmapsto -\rho(X) \\
Y & \longmapsto -\rho(Y) \end{cases}
\end{equation*}
The corresponding action on characters is:
\begin{align*}
(\sigma_{01})_*:\bmatrix \xi \\ \eta \\ \zeta \endbmatrix
& \longmapsto \bmatrix \xi \\ -\eta \\ -\zeta \endbmatrix \\
(\sigma_{10})_*:\bmatrix \xi \\ \eta \\ \zeta \endbmatrix
& \longmapsto \bmatrix -\xi \\ \eta \\ -\zeta \endbmatrix \\
(\sigma_{11})_*:\bmatrix \xi \\ \eta \\ \zeta \endbmatrix
& \longmapsto \bmatrix -\xi \\ -\eta \\ \zeta \endbmatrix. 
\end{align*}
Note that the sign change automorphisms are not induced by
automorphisms of $\pi$. 
\par
Thus the group $\Aut(\C^3, \kappa)$ is isomorphic to a semidirect
product 
\begin{equation*}
\pglz \ltimes (\Z/2\oplus\Z/2)
\end{equation*}
See Goldman~\cite{Gold:pTorus}, \S1.3 and Appendix, for a detailed
discussion. 
\par
The action of $\Aut(\C^3, \kappa)$ on $\C^3$ restricts to an action on
the $\islr$-character variety. However the subgroup $\mathfrak{S}_3$
does not preserve the individual components $\cv{j}{k}$. For example,
the automorphism induced by transposition of $x$ and $y$ maps
characters in $\cv{1}{0}$ to  characters in $\cv{0}{1}$, and
stabilizes $\cv{1}{1}$. For each component $\cv{j}{k}$ there is a
subgroup of finite index of $\Aut(\C^3, \kappa)$ that preserves
$\cv{j}{k}$. We call this subgroup \emph{the modular group of
$\cv{j}{k}$} and denote it by $\Gamma_{jk}$. 
\par
For instance, $\Gamma_{11}$ is generated by the quadratic reflections
\begin{equation*}
Q_x:\,
\left[\begin{matrix}
x \\ y \\ z
\end{matrix} \right]
\longrightarrow
\left[\begin{matrix}
yz -x \\ y \\ z
\end{matrix} \right],\quad
Q_y:\,
\left[\begin{matrix}
x \\ y \\ z
\end{matrix} \right]
\longrightarrow
\left[\begin{matrix}
x \\ xz - y \\ z
\end{matrix} \right],\quad
Q_z:\,
\left[\begin{matrix}
x \\ y \\ z
\end{matrix} \right]
\longrightarrow
\left[\begin{matrix}
x \\ y \\ -xy - z
\end{matrix} \right],
\end{equation*}
the sign-change automorphisms $\sigma_{jk}$, and the transposition
\begin{equation*}
t_{xy}:\,
\left[\begin{matrix}
x \\ y \\ z
\end{matrix} \right]
\longrightarrow
\left[\begin{matrix}
y \\ x\\ z
\end{matrix} \right]
\end{equation*}
where $x,y,z\in\R$. Since 
$\Gamma_{01} \cong \Gamma_{10} \cong \Gamma_{11}$,  
whenever it is clear from the context, we will drop the 
subscript and use $\Gamma$ to denote the relevant modular group
for the given component.
\subsection{Goldman's Result on the Real $\mop{SL}(2)$-Characters of the
Punctured Torus}
Recently, Goldman~\cite{Gold:pTorus} studied the action of the modular
group on the $\slr$ component of $\E{X}$. In this case
the modular group $\Gamma$ is isomorphic to 
\[
\pglz \ltimes (\Z/2\oplus\Z/2)
\]
We summarize Goldman's results in the theorem below.
\begin{thm}[Goldman]
Let $\kappa(\xi,\eta, \zeta) = \xi^2 + \eta^2 + \zeta^2 - \xi\eta\zeta
-2$ and let $c\in\R$. 
\begin{itemize}
\item For $c<-2$, the group $\Gamma$ acts properly and freely on
$\kappa\inv(c)\cap\R^3$;
\item For $-2 \leq c < 2$, there is a compact connected component $K_c$
of $\kappa\inv(c)\cap\R^3$ and $\Gamma$ acts properly and freely on
the complement $\kappa\inv(c)\cap\R^3 - K_c$;
\item For $18 \geq c > 2$, the group $\Gamma$ acts ergodically on 
$\kappa\inv(c)\cap\R^3$;
\item For $c> 18$, the group $\Gamma$ 
acts properly and freely on an open subset 
$\Omega_c\subset \kappa\inv(c)\cap\R^3 $ 
and acts ergodically on the complement of $\Omega_c$.
\end{itemize}

\end{thm}

\section{Fricke Spaces Within the Character Variety}\label{S:fricke}
We now focus on the components $\cv{j}{k}$ of the $\islr$-character
variety of $\pi$ that correspond to non-zero Steifel-Whitney classes
in $\Z_2\times\Z_2$. Our goal is to identify the regions in
$\cv{j}{k}$ that parametrize the Fricke spaces of the two
non-orientable surfaces whose  fundamental group is isomorphic to
$\pi$, namely \PKB and \PMB (see section~\ref{SS:surf-rank-two}). 
For short, we call these surfaces the \emph{PK-bottle} and the
\emph{PM-band} respectively.
\subsection{The Fricke Space of the Punctured Klein-Bottle}
We first derive a presentation of the fundamental group of the
PK-bottle $\PKB$ from a presentation of $\pi_1(\pants)$ via the 
Higman-Neumann-Neumann (HNN) extension construction. 
Let $\delta$, $\alpha'$, and $\alpha''$
be the homotopy classes in $\pi_1(\pants)$ of the boundary loops of $\pants$.
Then $\pi_1(\pants)$ admits the following (redundant) presentation:
\begin{equation*}
\pi_1(\pants) = \langle \delta, \alpha', \alpha'' \mid 
\delta \alpha' \alpha'' = I \rangle
\end{equation*}
The PK-bottle $\PKB$ is obtained from the pair of pants $\pants$ via
identifying the boundary components $\alpha'$ and $\alpha''$ by a (orientation
reversing) diffeomorphism $\phi$ (see Fig~\ref{F:hnnpkb}). Let $H = \langle \alpha' \rangle$ be the
cyclic subgroup of $\pi_1(\pants)$ generated by $\alpha'$. Then the induced
homomorphism $\phi_*$ sends $\alpha'$ to $\alpha''$ and thus its
restriction to $H$ is a monomorphism:
\[ 
\phi_*|_H : H \longmapsto \pi_1(\pants)
\]
The corresponding HNN-extension of $\pi_1(\pants)$ is given by:
\begin{equation}
\pi_1(\pants) *_{\phi} H = \langle \delta, \alpha', \alpha'', \beta \mid \delta \alpha' \alpha'' = I, 
\phi_*(\alpha') := \alpha'' = \beta \alpha' \beta^{-1} \rangle
\end{equation}
This defines a presentation of $\pi_1(\PKB)$ which reduces to the more
familiar: 
\begin{equation}
\pi_1(\PKB) = \langle \alpha, \beta, \delta \mid 
\alpha\beta\alpha \beta^{-1} = \delta\inv \rangle
\end{equation}
\begin{figure}[htbp]{\label{F:hnnpkb}}
\begin{center}
\input{hnn-pkbot.pstex_t}
\end{center}
\caption{\PKB\ as a quotient space of \pants}
\end{figure}
Both $\pi_1(\pants)$ and $\pi_1(\PKB)$ are free of rank two. In particular, 
$\pi_1(\pants)$ is freely generated by $\alpha'$ and $\alpha''$, 
while $\pi_1(\PKB)$
is freely generated by $\alpha$ and $\beta$. The quotient map 
\begin{equation*}
\pi_1(\pants) \lrarw \pi_1(\pants) *_\phi H 
\end{equation*}
defines a monomorphism of
the fundamental groups:
\begin{align*}
\phi_*:\pi_1(\pants) & \lrarw \pi_1(\PKB) \\
\alpha' &\longmapsto \alpha \\
\alpha''& \longmapsto \beta \alpha \beta^{-1} \\
\delta &\longmapsto \beta\alpha^{-1}\beta^{-1}\alpha^{-1}
\end{align*}
Let $\rho \in \E{R}_{0,1}$ and let $X\in G_0$, $Y\in G_1$ denote the
images under $\rho$ of $\alpha$  and $\beta$ respectively. Then the
images under $\rho\circ\phi_*$ of $\alpha'$ and $\alpha''$ will be
$X$ and $YXY^{-1}$ respectively. Thus the induced map $\phi^*$ of
representation spaces is:  
\begin{align*}
\phi^*: \E{R}_{0,1} \cong  G_0 \times G_1 & \lrarw 
G_0 \times G_0 \cong \E{R}_{0,0} \\
(X, Y) & \longmapsto (X, Y X Y^{-1})
\end{align*}
A representation of $\pi_1(\pants)$ defined by $(X,Z) \in G_0 \times G_0
\cong \rep{0}{0}$
pulls back to a representation of $\pi_1(\PKB)$ in $\rep{0}{1}$ if and only if $Z$ is in
the centralizer of $X$ in $G_1$, that is, there exists $Y \in G_1$, such
that $Z = YXY^{-1}$. 
The map $\phi^*$ descends to a map (denoted again, by abuse of
notation, as $\phi^*$) between the corresponding components of the 
character varieties:
\begin{align*}
\phi^*: \E{X}_{0,1} \cong \R\times i\R \times i\R  & \lrarw 
\R^3  \cong \E{X}_{0,0} \\
\left[\begin{matrix} 
\xi \\ \eta \\ \zeta
\end{matrix} \right]
& \longmapsto 
\left[\begin{matrix} 
\xi \\ \xi \\ -\eta^2 -\zeta^2 + \xi\eta\zeta + 2
\end{matrix} \right]
\end{align*}
Setting $\xi = x$, $\eta = iy$, and $\zeta = iz$ with $x,y,z \in \R$
establishes an isomorphism \mbox{$\R\times i\R\times i\R \cong \R^3$} such that
the map $\phi^*$ can be expressed as follows:
\begin{align*}
\phi^*: \E{X}_{0,1} \cong \R^3  & \lrarw \R^3 \cong \E{X}_{0,0} \\
\left[\begin{matrix} 
x \\ y \\ z
\end{matrix} \right]
& \longmapsto 
\left[\begin{matrix} 
x \\ x\\  y^2 +z^2 -xyz + 2
\end{matrix} \right]
\end{align*}
\nl
A representation $\rho \in \E{R}_{0,0}$ gives rise to a hyperbolilc
structure on $int(\pants)$  if and only if $\rho(\alpha')$,
$\rho(\alpha'')$ and $\rho(\delta)$ are hyperbolic with
non-intersecting axes that span an
ultra-ideal triangle in $\Ht$. This condition is equivalent to 
$\chi([\rho]) \in (-\infty, -2)^3 \cup (2, +\infty)^3$, such that odd
number of components of $\chi([\rho])$ are negative, or equivalently
\begin{align}
|x| &> 2\label{E:xgt} \\
y^2 +z^2 - xyz +2 &< -2 \label{E:deltatr}
\end{align}
But characters in $\cv{0}{1}$  that lie in the fibers of $\phi^*$ 
correspond to discrete embeddings of $\pi_1(\PKB)$ in $\rep{0}{1}$ with
quotient a PK-Bottle, and thus to marked
hyperbolic structures on $int(\PKB)$. Note that inequality~\ref{E:deltatr}
already implies that $|x|>2$. We have thus proved the first part of
the following
\begin{prop}\label{P:Kbottle}
Let $x, y, z\in \R$ be the coordinate functions on $\E{X}_{0,1}$. 
\begin{enumerate}
\item The region $\Omega_0^K \subset \R^3$ parametrized by the inequality
\begin{equation*} 
y^2 + z^2 - xyz + 4 < 0 
\end{equation*}
corresponds to discrete geometric $K$-embeddings inside 
\rep{0}{1}. 
\item 
Let $\Gamma_K$ be the image of the mapping class group
$\mop{Map}({\PKB}))$ in $\Gamma$ under the 
Nielsen homomorphism and let $\Gamma_0^K = \Gamma/\Gamma_K$
be the coset space of $\Gamma_K$ in $\Gamma$.
The space of all discrete $K$-embeddings identifies with 
\[
\Omega^K = \coprod_{\gamma \in \Gamma_0^K}\gamma\Omega_0^K
\]
\end{enumerate}
\begin{proof}
To prove part (2) we must show that if for some $\gamma\in\Gamma$ the
intersection $\Omega_0^K\cap \gamma\Omega_0^K$ is nonempty, then
$\gamma\in \mop{Map}(\PKB)$.Suppose that $[\rho]\in\Omega_0^K\cap 
\gamma\Omega_0^K$. The automorphism $\gamma$ of the character space
$\C^3$ is induced by an automorphism $\tilde{\gamma}$ of
$\pi_1(\PKB)$ such that 
\[
[\rho_0\circ\tilde{\gamma}] = \gamma\circ[\rho]
\]
for some $\rho_0 \in \Omega_0^K$. Thus $\rho_0\circ\tilde{\gamma}$ is
also a discrete geometric $K$-embedding with respect to the original
peripheral structure. This means that $\tilde{\gamma}$ preserves the
conjugacy class of each boundary component, and in particular that 
it represents an equivalence class in the mapping class group of
$\PMB$. 
\end{proof}
\end{prop}
Once we have computed the Fricke space of the PK-bottle inside
$\cv{0}{1}$, it is easy to obtain the Fricke spaces sitting  in the
other two non-zero components of the $\islr$-character variety.
We show how to do this for $\cv{1}{1}$. The class of the automorphism 
\[
\psi:\,\alpha \longmapsto \alpha\beta\inv
\]
of $\pi$ induces the transposition
\[
t_{xz}:\,
\left[\begin{matrix}
x \\ y \\ z
\end{matrix} \right]
\longrightarrow
\left[\begin{matrix}
z \\ y\\ x
\end{matrix} \right]
\]
on the character variety $\E{X}$. In particular $t_{xz}$ interchanges
$\cv{0}{1}$ and $\cv{1}{1}$ and consequently the image of
$\Omega_0^K$ in $\cv{1}{1}$ will be given by
\begin{equation}
x^2 + y^2 - xyz + 4 < 0
\end{equation}
This corresponds to a change of marking, with respect to which
the new presesntation of $\pi_1(\PKB)$ becomes 
\[
\langle \gamma, \beta, \delta \mid \gamma\beta^2\gamma = 
\delta\inv \rangle
\]
\subsection{The Fricke Space of the Punctured M\"obius Band}
\begin{figure}[htbp]{\label{F:qp-sphere}}
\begin{center}
\input{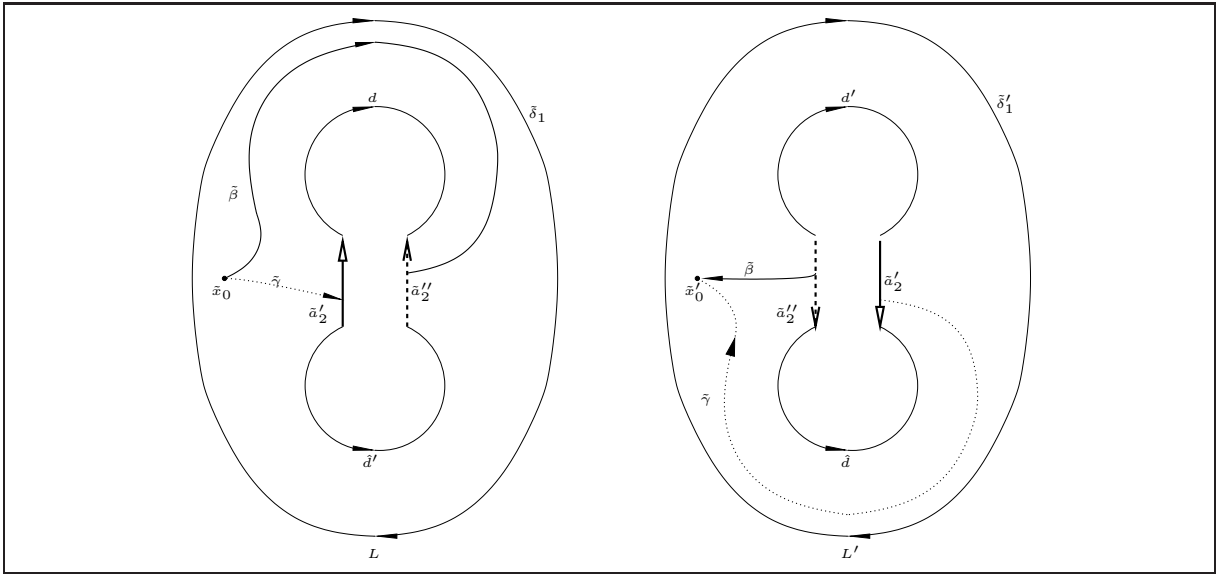}
\end{center}
\caption{The Two Sheets of \qps{} as a Double Cover of \PMB}
\end{figure}

We now determine the characters of the discrete geometric embeddings
inside the $(1,1)$-component of \linebreak[1]
$\hom{\pi}{\islr}$. To this end,  consider the
inclusion of the fundamental group of the quadruply punctured sphere
$\pi_1(S_{0,4})$ into $\pi_1(\PMB)$ induced by the (double) covering
map 
\[
q: S_{0,4} \lrarw \PMB
\]
Suppose $\pi_1(\PMB)$ is given the presentation discussed in
section (\ref{SS:surf-rank-two}), namely:
\begin{equation}
	\pi = \langle \beta, \gamma, \delta_1, \delta_2 \mid \delta_1 =
	\beta \gamma, \delta_2 = \beta\gamma\inv \rangle
\end{equation}
The orientable double cover of $\PMB$, which we denote for the
moment by $\PMBt$, is a connected surface
of genus $0$ with $4$ boundary components and thus homeomorphic to 
$\qps$. The two
sheets of the covering map, denoted $L$ and $L'$, are represented
schematically in Fig.~\ref{F:qp-sphere}.
To obtain the covering surface, $L$ and
$L'$ are glued via identifying the
two pairs of edges, respectively marked with $\tilde{a}_2'$ and 
$\tilde{a}_2''$ (these are respectively the lifts of $a_2'$ and
$a_2''$ from Fig. \ref{F:pM-band}). After this identification, the geodesic segments 
$d$ and $\hat{d}$ concatenate to form a boundary loop, denoted by
$\t{\delta}_2$. Similarly, the geodesic segments $d'$ and
$\hat{d}'$ concatenate to form a boundary loop, denoted by
$\t{\delta}_2'$. Thus the boundary of $\PMBt$  consists of four 
components, namely $\t{\delta}_1$, $\t{\delta}_1'$,
$\t{\delta}_2$, and $\t{\delta}_2'$. Their corresponding 
elements in the fundamental group of $\PMBt$, denoted again by
$\t{\delta}_1$, $\t{\delta}_1'$, $\t{\delta}_2$, and
$\t{\delta}_2'$, are in fact the lifts of $\delta_1,\,\delta_2 \in
\pi_1(\PMB)$. 
\\*[1.5\baselineskip]
More precisely, fix a base point $x_0 \in \PMB$ and let $\tilde{x}_0$
be the lift of $x_0$ to $L$. Consider the fundamental group
$\tilde{\pi}$ of $\PMBt$ with base point $\tilde{x}_0$. Let
$\tilde{x}_0'$ be the lift of $x_0$ to $L'$. In $\PMBt$ there are two
homotopy classes of paths with initial point $\tilde{x}_0$
and end point $\tilde{x}_0'$; one class projects to the generator
$\beta \in \pi$, the other one projects to $\beta\inv$. Let
$\t{\beta}$ and $\t{\beta}'$ be representatives of each
class respectively. Similarly, let $\t{\gamma}$ and
$\t{\gamma}'$ be representatives of the two homotopy classes of
paths that project to $\gamma$ and $\gamma\inv$ respectively (see
Fig. (\ref{F:qp-sphere})). With this notation, the homotopy classes of
the boundary components of $\PMBt$ can be written as follows:
\begin{equation*}
\begin{aligned}
\t{\delta}_1 &= [\t{\beta}*(\t{\gamma}')\inv] \\
\t{\delta}_2 &= [\t{\beta}*\t{\gamma}\inv] \\
\end{aligned}
\quad\quad
\begin{aligned}
\t{\delta}_1' &= [\t{\gamma}*(\t{\beta}')\inv] \\
\t{\delta}_2' &= [\t{\gamma}'*(\t{\beta}')\inv]
\end{aligned}
\end{equation*}
where '$*$' denotes path concatenation.
From these equations, the images of $\t{\delta}_1, \t{\delta}_1',
\t{\delta}_2$, and $\t{\delta}_2'$ under the induced covering
projection map $q_*$ can be easily computed:
\begin{equation*}
\begin{aligned}
q_*: \pi_1(\PMBt) &\lrarw \pi_1(\PMB)  \\
\left[
\begin{matrix}
\t{\delta}_1  \\
\t{\delta}_2  \\
\t{\delta}_1'  \\
\t{\delta}_2' \\
\end{matrix}
\right]
& \longmapsto
\left[
\begin{matrix}
\beta\gamma \\
\beta\gamma\inv \\
\gamma\beta \\
\gamma\inv\beta 
\end{matrix}
\right]
\end{aligned}
\end{equation*}
\nl
But $\t{\delta}_1, \t{\delta}_1', \t{\delta}_2$, and $\t{\delta}_2'$
satisfy the relation
\[
\t{\delta}_1 \t{\delta}_2' (\t{\delta}_1')\inv  (\t{\delta}_2)\inv = I
\]
and thus if we set 
\[
A = \t{\delta}_1, \quad B = \t{\delta}_2', \quad C= (\t{\delta}_1')\inv,
\quad  D = (\t{\delta}_2)\inv
\]
we obtain a presentation of $\pi_1(\PMBt)$ such that 
\begin{equation*}
\pi_1(\PMBt) \equiv \pi_1(S_{0,4}) = \langle A,B,C,D \mid ABCD = I \rangle
\end{equation*}
The induced map $q_*$ is then uniquely determined by the images of the three
generators $A$, $B$, and $C$
\begin{equation}\label{E:qps-pmb}
q_*:
\left[
\begin{matrix} 
A \\ B \\ C 
\end{matrix}
\right]
\longmapsto
\left[
\begin{matrix}
\beta\gamma \\ \gamma^{-1}\beta \\ \beta^{-1}\gamma^{-1}
\end{matrix}
\right]
\end{equation}
The image of $q_*$ in $\pi_1(\PMB)$ is the index-$2$ subgroup 
consisting of words of even length, and generated by $\beta\gamma$, $\gamma\inv \beta$,
and $\beta\inv \gamma\inv$. The elements of $\pi_1(\PMBt)$ whose traces generate
the coordinate ring of the character variety of \PMBt, are (see for
instance Magnus~\cite{Magnus})
\[
A, B, C, D, AB, BC, CA
\]
From (\ref{E:qps-pmb}) we compute the images under $q_*$ of
the last four elements:
\begin{equation*}
q_*:
\left[
\begin{matrix} 
D \\ AB \\ BC \\ CA
\end{matrix}
\right]
\longmapsto
\left[
\begin{matrix}
\beta\gamma\inv \\ \beta^2 \\ \gamma^{-2} \\ \beta\inv \gamma\inv \beta\gamma
\end{matrix}
\right]
\end{equation*}
Thus the induced map $q^*$ on representation spaces:
\begin{align*}
q^*: \hom{\pi_1(\PMB)}{G} & \lrarw \hom{\pi_1(\qps)}{G} \\
	\rho & \longmapsto \rho \circ q_*
\end{align*}
can be expressed in terms of the images $X = \rho(\beta)$ and 
$Y = \rho(\gamma)$ of the generators $\beta$ and $\gamma$ 
as follows
\begin{align*}
q^*: \hom{\pi_1(\PMB)}{G} & \lrarw \hom{\pi_1(\qps)}{G} \\
\left[
\begin{matrix}
X \\ Y
\end{matrix}
\right] 
& \longmapsto
\left[
\begin{matrix}
XY \\ Y\inv X \\ X\inv Y\inv \\ XY\inv \\ X^2 \\ Y^{-2} \\ X\inv Y\inv XY
\end{matrix}
\right]
\end{align*}
Consequently, the induced map on the $\slc$-character varieties, denoted
again by $q^*$, is defined in terms of trace coordinates as follows:
\begin{equation}\label{E:q-star-pMband}
\begin{aligned}
q^*: \hom{\pi_1(\PMB)}{G}\doubleslash G & \lrarw \hom{\pi_1(\qps)}{G}\doubleslash G \\
\C^3 \ni
\left[
\begin{matrix}
\xi \\ \eta \\ \zeta
\end{matrix}
\right] 
& \longmapsto
\left[
\begin{matrix}
\zeta \\ \xi\eta - \zeta\\ \zeta \\ \xi\eta - \zeta\\ \xi^2 -2  \\ 
\eta^2 -2  \\ \kappa(\xi,\eta, \zeta)
\end{matrix}
\right] \in \C^7
\end{aligned}
\end{equation}
where 
\[
\xi = \mop{tr}(X), \; \eta = \mop{tr}(Y), \; \zeta = \mop{tr}(XY)
\]
Recall that the $(1,1)$-component $\cv{1}{1}$ of the $\islr$-character
variety of $\PMBt$ is isomorphic to $i\R\times i\R \times \R$, and
therefore identifies with $\R^3$ by setting
\[
\xi = ix, \quad \eta = iy, \quad \zeta = z
\]
with $x,y,z \in \R$. The restriction of $q^*$ to $\cv{1}{1}$ 
can be expressed in terms of the $x,y,z$ coordinates in the following
way:
\begin{equation}\label{E:qstar-pmband}
\begin{array}{rrcc}
q^*:& \cv{1}{1}\cong \R^3  & \lrarw & \R^7 \\
    &
\left[
\begin{matrix}
x \\ y \\ z
\end{matrix}
\right] 
& \longmapsto &
\left[
\begin{matrix}
z \\ - xy - z\\ z \\ -xy - z\\ -x^2 -2  \\ -y^2 -2  \\ \kappa(ix,iy, z)
\end{matrix}
\right] 
\end{array}
\end{equation}
The coordinate ring of the $\slc$-character variety of $\qps$ is
generated by the traces
\[ 
\begin{gathered}
a = \tr({A}),\; b = \tr({B}),\;  c=\tr({C}), \;d = \tr({D}) \\
t_{AB} = \tr({AB}), \; t_{BC} = \tr({BC}), \; t_{CA} = \tr({CA})
\end{gathered}
\]
subject to the relation (see for instance Magnus (\cite{Magnus})):
\begin{align*}
t_{AB}^2 + t_{BC}^2 + t_{CA}^2 & + t_{AB}t_{BC}t_{CA}  \\
	& = (ab+cd)t_{AB} + (ad +bc)t_{BC} + (ac + bd)t_{CA} \\
	    & \qquad - (a^2 + b^2 + c^2 + d^2 + abcd - 4).
\end{align*}
Thus the map
\begin{equation*}
\E{X}_{\qps}^{\C} = \hom{\pi_1(\qps)}{\slc}\doubleslash \slc\longrightarrow \C^7
\end{equation*}
defined by $(a, b, c, d, t_{AB}, t_{BC}, t_{CA})$ embeds $\E{X}$ onto
a hypersurface in $\C^7$. The set of real points of
$\E{X}_{\qps}^{\C}$ parametrizes equivalence classes of
representations in $\slr$. The following Proposition (compare
Keen~\cite{Keen}, Theorem 2) provides necessary and sufficient
conditions for
\begin{equation*}
(a,b,c,d,t_{AB}, t_{BC})\in\R^6 
\end{equation*}
to represent an equivalence
class of discrete geometric $Q$-embeddings of $\pi_1(\qps)$ into $\slr$. 
\begin{prop}
Given $a< -2$, $b < -2$, $c < -2$, $d<-2$, $t_{AB}< -2$, and $t_{BC} <
-2$,  there exist elements $A^*,B^*,C^*, D^* \in \slr$ such that 
\begin{equation*}
\begin{gathered}
a = \tr({A^*}),\; b = \tr({B^*}),\;  c=\tr(C^*), \;d = \tr(D^*) \\ 
t_{A^*B^*} = \tr(A^*B^*), \; t_{B^*C^*} = \tr(B^*C^*)
\end{gathered}
\end{equation*}
and such that $F = \langle A^*, B^* \rangle$ and 
$F' = \langle C^*, D^*\rangle$ are Fuchsian groups of signature
$(0;3)$. Moreover, the group $H = F*_{\Z} F'$, which is the
amalgamated product of $F$ 
and $F'$ over the cyclic subgroups 
\begin{equation*}
 \Z\cong\langle A^*B^* \rangle \subset F, \qquad \Z\cong\langle C^*D^* \rangle
\subset F',
\end{equation*}
is Fuchsian and represents a \emph{marked} surface of signature
$(0;4)$. Every such marked surface is so representable.
\end{prop}
But a Fuchsian group $H$, whose quotient $\Ht/H$ is a
marked hyperbolic surface $S$, gives rise to a discrete geometric
embedding $\rho \in \hom{\pi_1(S)}{\slr}$ determined (up to
conjugation) by the correspondence between the generators of
$\pi_1(\qps)$ and those of $H$. Thus we have the following
\begin{cor}\label{C:Keen-Me}
Given $a, b,c,d,t_{AB}, t_{BC} \in \R$, such that 
\begin{equation}\label{E:Keen-ineq}
a< -2,\; b < -2,\; c < -2,\; d<-2,\; t_{AB}< -2,\; t_{BC} < -2, 
\end{equation} 
there exists a hypebolic surface $\qps$ of signature $(0;4)$ with
marking $\Lambda =\{A,B,C,D\}$ and a discrete $Q$-geometric embedding
$\rho \in \hom{\pi}{\slr}$, such that 
\begin{equation*}
\begin{gathered}
a = \tr({\rho(A)}),\; b = \tr({\rho(B)}),\;  c=\tr({\rho(C)}), \;d =
\tr({\rho(D)}) \\ 
t_{AB} = \tr(\rho(AB)), \; t_{BC} = \tr(\rho(BC))
\end{gathered}
\end{equation*}
All such discrete $Q$-geometric embeddings are obtained this way.
\end{cor}
It follows from Corollary~\ref{C:Keen-Me} that a real point on the
character variety $\E{X}_{\qps}^{\C}$, that is,. a point $p =
(a,b,c,d,t_{AB}, t_{BC}, t_{CA}) \in \R^7$ whose coordinates satisfy
the Fricke relation, corresponds to a discrete geometric $Q$-embedding
if and only if its coordinates satisfy
inequalities~(\ref{E:Keen-ineq}). 
\par
This conclusion together with the observation that some of the
expressions that define the image of $q^*$ in (\ref{E:qstar-pmband}) 
satisfy the given inequalities trivially, prove the first part of the
following 
\begin{prop}
Let $x, y, z\in \R$ be the coordinate functions on $\E{X}_{1,1}$. 
\begin{enumerate}
\item The region $\Omega_0^M$ of $\R^3$ parametrized by the
inequalities: 
\begin{align*} 
xy + z &> 2 \\
z & < - 2 
\end{align*}
corresponds to discrete geometric $M$-embeddings  inside 
\rep{1}{1}. 
\item 
Let $\Gamma_M$ be the image of the mapping class group
$\mop{Map}({\PMB}))$ in $\Gamma$ under the 
Nielsen homomorphism and let $\Gamma_0^M = \Gamma/\Gamma_M$
be the coset space of $\Gamma_M$ in $\Gamma$.
The space of all discrete $M$-embeddings identifies with 
\[
\Omega^M = \coprod_{\gamma \in \Gamma_0^M}\gamma\Omega_0^M
\]
\end{enumerate}
\end{prop}
The proof of the second part of the Proposition is analogous to the
proof of Proposition~\ref{P:Kbottle}

\section{Classification of Characters and the $\tau$-reduction Algorithm}

The next two sections discuss the $\Gamma$-action 
on the components 
\begin{equation*}
\cv{0}{1} \cup \cv{1}{0} \cup \cv{1}{1}
\end{equation*}
of the $\islr$-character variety of $\pi$.  
We shall mainly work in $\cv{1}{1}$, but the results carry
over verbatim to $\cv{0}{1}$ and $\cv{1}{0}$.  
\par
\renewcommand{\baselinestretch}{1.5} 

Denote by ``$\sim$'' the equivalence relation induced by the 
$\Gamma$-action on characters: 
For $u,v \in \cv{1}{1}$
\begin{equation*}
u\sim v \Longleftrightarrow \exists
\gamma \in \Gamma \text{~ such that~} \gamma u = v
\end{equation*}
Since $\kappa$ is $\Gamma$-invariant, 
$u\sim v$ implies $\kappa(u) = \kappa(v)$. 
In this section we show that when $c<-14$ there exist essentially two 
types of equivalence classes of characters upon each of which $\Gamma$ 
acts in a substantially different manner.

\begin{thm}\label{T:classify}
Suppose that $u \in \cv{1}{1}$ satisfies $\kappa(u) < -14$. 
Then $\exists u'=(x',y',z') \sim u$ 
such that either:
\begin{description}
\item[M] $u' \in \oM$, in which case $u'$ (and therefore $u$) is a character
of a Fuchsian representation whose quotient is homeomorphic to a 
a once-punctured M\"obius band; or 
\item[E] $\zbar' \in (-2, 2)$ in which case $u'$ is the character of a
representation mapping the peripheral element $\beta\gamma\inv\in\pi$ to an
elliptic element of $\slr$.
\end{description}
\end{thm}

\subsection{Notation}
For any $(x,y,z) \in \R^3$, let 
\begin{equation*}
\zbar := -xy -z 
\end{equation*}
Since the quadratic reflection
\[
Q_z:\;
\left[\begin{matrix}
x \\ y \\ z
\end{matrix} \right]
\longrightarrow
\left[\begin{matrix}
x \\ y \\ -xy - z
\end{matrix} \right]
\]
preserves every $\kappa\inv(c)$ for $c \in \R$, it interchanges 
$z$ and $\zbar$. Fix $x, y \in \R$. 
Then $z$ and $\zbar$
are the two (necessarily real) roots of  the quadratic polynomial
\begin{equation}\label{E:zquad}
z^2 + (xy)z - x^2 - y^2 - 2 - c
\end{equation}
since
\begin{equation}\label{E:viet}
\begin{aligned}
z + \zbar & = -xy \\
z\zbar &= -x^2 - y^2 -2 -c.
\end{aligned}
\end{equation}
Thus $Q_z$ is the deck transformation 
of the double covering 
\begin{equation*}
\Pi:\kappa\inv(c)\longrightarrow \Pi(\kappa\inv(c))
\end{equation*}
obtained by restriction of the 
projection $\Pi(x,y,z)= (x,y)$ to the $xy$-plane.

Denote by $\ec$ the region 
\[
\{(x,y,z)\in \kappa\inv(c)\mid \zbar \in (-2, 2)\}
\]

\subsection{The level sets of $\kappa$ and the Fricke space of $\PMB$}
Next we find a necessary condition for \oM and $\kappa\inv(c)$ to have
a non-empty intersection. 
\begin{lemma}\label{L:omegaM}
Suppose $u=(x,y,z) \in \oM \cap \kappa\inv(c)$. Then $c<-14$
\begin{proof}
Since $u\in\oM$, by defintion $z < -2$, and $\zbar < -2$. Therefore
$z + \zbar < -4$ and $z\zbar > 4$. Then (\ref{E:viet}) imply that
\begin{equation}\label{E:xyineq}
\begin{aligned}
xy & > 4 \\ x^2 + y^2 &< -c - 6
\end{aligned}
\end{equation}
Since any $x,y \in \R$ satisfy $x^2 + y^2 \ge 2xy$,
\begin{equation*}
c = -6 - (-c - 6) \le -6 - (x^2 + y^2) \le -6 - 2xy < -14. 
\end{equation*}
\end{proof}
\end{lemma}
\begin{cor}
$\oM \cap \kappa\inv(c) \neq \varnothing$ if and only if $c < -14$.
\end{cor}
\begin{remark} The sign change transformations
\begin{equation}
\sigma_{xz}: \left[
\begin{matrix} x\\ y\\ z \end{matrix} 
\right] \longmapsto
\left[
\begin{matrix} -x \\ y \\-z \end{matrix}
\right],\quad 
\sigma_{yz}: \left[
\begin{matrix} x\\ y\\ z \end{matrix} 
\right] \longmapsto
\left[
\begin{matrix} x \\ -y \\-z \end{matrix}
\right]
\end{equation}
preserve $\kappa$ and commute with $Q_z$. Since
\begin{equation*}
\zbar(\sigma_{xz}(u)) =\zbar(\sigma_{yz}(u)) = -\zbar(u) 
\end{equation*}
it suffices to consider
$u \in \cv{1}{1}$ such that $\kappa(u) < -14$, $z < -2$ and $\zbar >2$. 
\end{remark}
The proof of Theorem~\ref{T:classify} proceeds in two steps: first we
construct a function $\tau$ 
that is non-decreasing along certain subsets of the orbit of each
character in \oM. Then, (always assuming  $c < -14$ and 
$u \notin \oM\cup\ec $), we find a finite sequence of characters
$u_0=u, u_1, 
\dots u_N$, and a constant $T>0$ such that $\tau(u_i) > \tau(u_{i+2}) +
T$, and either $u_N \in \oM$, or $u_N \in \ec$. 

\subsection{The Quadratic Reflections}
Recall that, apart from $Q_z$, two other quadratic reflections act on
characters as automorphisms of $\kappa$:
\begin{equation*}
Q_x:\;
\left[\begin{matrix}
x \\ y \\ z
\end{matrix} \right]
\longrightarrow
\left[\begin{matrix}
yz -x \\ y \\ z
\end{matrix} \right] \quad \text{and} \quad
Q_y:\;
\left[\begin{matrix}
x \\ y \\ z
\end{matrix} \right]
\longrightarrow
\left[\begin{matrix}
x \\ xz - y \\ z
\end{matrix} \right]
\end{equation*}
They are induced, respectively by the following automorphisms of $\pi$
(cf.~Goldman~\cite{Gold:pTorus})
\begin{equation*}
q_x:\;
\begin{array}{rll}
X & \longmapsto & (XY)Y^2X(Y\inv X\inv) \\ 
Y & \longmapsto & (XYX\inv)Y\inv(XY\inv X\inv)\\ 
XY & \longmapsto & (XY^2)(XY)(Y^{-2}X\inv) 
\end{array}
\end{equation*}
and
\begin{equation*}
q_y:\;
\begin{array}{rll}
X & \longmapsto &(XY)X(Y\inv X\inv) \\ 
Y & \longmapsto &(XY)X\inv Y\inv X\inv(Y\inv X\inv)  \\ 
XY & \longmapsto & (XY)\inv 
\end{array} 
\end{equation*}
On the other hand $Q_z$ is induced by the automorphism
\begin{equation*}
q_z:\;
\begin{array}{rll}
X & \longmapsto & X \\ 
Y & \longmapsto & Y\inv \\
XY & \longmapsto & XY\inv 
\end{array} 
\end{equation*}
which interchanges the two boundary components of $\PMB$. Therefore
its coset in $\Out(\pi)$ lies in the image of the mapping class group
under the Nielsen homomorphism
\eqref{Eq:Nielsen}. In contrast $q_x$ and $q_y$ do not
represent elements of the mapping class group of $\PMB$. 
\par
Clearly $Q_z$ preserves $Q_x$ and 
$Q_y$ preserves \oM.
Indeed, for each $u\in \oM$
\[
\zbar(Q_x(u)) = -y(yz -x) -z = -z(y^2+1) + xy >6
\]
since (\ref{E:xyineq}) implies $xy > 4$ and $z<-2$.
Similarly, 
\[\zbar(Q_y(u)) > 6
\]
and therefore $Q_x(u), Q_y(u) \notin \oM \cup \ec$.

\subsection{The orbit as a binary tree}
Consider the subgroup $\Lambda \subset \Gamma$ generated by 
$Q_x, Q_y, Q_z$.
Since $\Lambda\cong \Z_2*\Z_2 *\Z_2$, 
it can be naturally associated with a trivalent tree, 
$T_{\Lambda}(V,E)$, defined as follows. Each node $v\in V$ represents a
group element, that is, a reduced word on the generators. Nodes 
corresponding to words $w_1$ and $w_2$ respectively, are
linked by an edge $e \in E$ if there is a generator 
$\lambda$ such that $w_1 = \lambda w_2$. 
\par
For any $u\in \oM$, the tree $T_{\Lambda}(V,E)$ imparts a binary
forest structure on the orbit $\Lambda\cdot u$
\begin{defin}
Let $u\in\oM$. Define a binary tree $\blmd$ inductively as follows
\begin{itemize}
\item $u$ is the root of $\blmd$
\item $Q_x(u)$ and $Q_y(u)$ are, respectively, 
the left and the right descendent of $u$.
\item Suppose $v$ is an arbitrary node and let $\hat{v}$ be its
parent. Let $\lambda \in \Lambda$ be the generator such that $v = \lambda
\hat{v}$. Then the descendents of $v$ are $\lambda_1 v$, and
$\lambda_2 v$, where $\lambda_1,\lambda_2 \neq \lambda$.
\end{itemize}
\end{defin}
For every $u \in \oM$ there is a ``dual'' tree $\blm(Q_z(u))$
rooted at $Q_z(u)$. Since 
$q_x, Q_y, Q_z$ freely generate $\Lambda$, the 
orbit $\Lambda\cdot u$ is the disjoint union of two binary trees
\[
\blmd \amalg \blm(Q_z(u))
\]
\begin{figure}[htbp]
\centering
\input{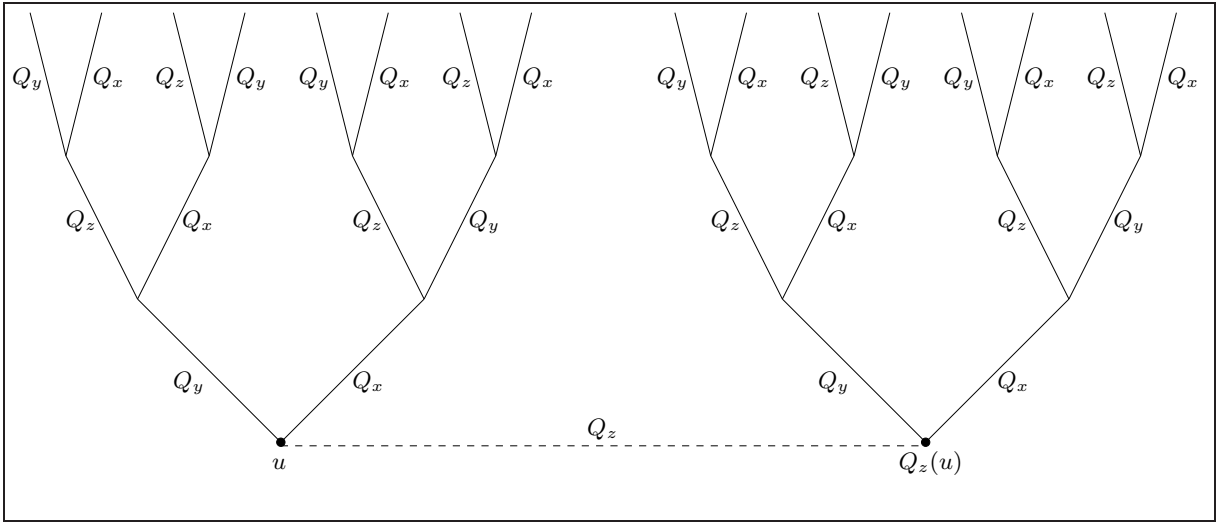}
\caption{The binary trees rooted at $u$ and $Q_z(u)$}
\end{figure}

\subsection{The $\tau$-function}
\begin{prop}\label{P:tautree}
For every $u\in \oM$, the function
\begin{align*}
\tau: \cv{1}{1} & \longrightarrow \R \\
\bmatrix x \\ y \\ z \endbmatrix & \longmapsto  -z\zbar
\end{align*}
does not decrease along the depth levels of \blmd. More precisely, if $v$ is
a node of \blmd, and $v_l, v_r$ are its left and right descendents
respectively, then 
\[ 
\tau(v) \leq \tau(v_l), \quad 
\tau(v) \leq \tau(v_r)
\]
with at least one of the inequalities strict. 
\end{prop}
\begin{proof}
Fix $u_0=(x_0, y_0, z_0)\in \oM$ and let $u_1 = (x_1, y_1, z_1) =
Q_x(u_0)$. By definition $z_0 < -2$ and $\zbar < -2$ and therefore
\[
\tau(u_0) = -z\zbar < 0
\]
Then, $x_0y_0 = -(z_0 + \zbar_0)$ implies 
\begin{equation}\label{E:tauroot}
\zbar_1 = -y_0^2 z_0 + x_0y_0 - z_0 = -z_0(y^2+2) - \zbar_0 > -2z_0 -
\zbar_0 > 0 
\end{equation}
and consequently:
\[
\tau(u_1) > 0 > \tau(u_0)
\]
Similar estimates apply verbatim in the case when $u_1 = Q_y(u_0)$ and
therefore $\tau$ is strictly increasing at the root of $\blm(u_0)$. 
However, the same type of algebraic argument does not extend to other
nodes of $\blm(u_0)$ directly. In order to
proceed with the induction step we analyse the orbits of $\Lambda$
geometrically. 
\par
Let $L_{x-x_0}, L_{y-y_0}, L_{z-z_0}$ 
denote the level sets of the $x$-, $y$-, and
$z$-coordinate functions respectively. For example
\[
L_{x-x_0} = \{(x,y,z) \in \R^3 \mid x - x_0=0\}
\]
Consider the conic $\hczz = L_{z-z_0} \cap \kappa\inv(c)$
and its projection $\Pi(\hczz)$ to the $xy$-plane:
\[
\Pi(\hczz) = \{(x,y)\in \R^2\mid k(x,y,z_0) = c\}
\]
\begin{lemma}\label{L:hypcz}
Let $c<2$ and $|z|>2$. Then $\hcz$ 
is a hyperbola with principal axes
\begin{equation}\label{E:principalaxes}
a_1:= \{(x,y) \mid  x+y = 0 \}, \quad
a_2:= \{(x,y) \mid  x-y = 0 \}
\end{equation}
and asymptotes
\begin{align}\label{E:equasmp}
l_1 & := \{(x,y) \mid  
y = \frac{z - \sqrt{z^2 -4}}{2} x \} \\
l_2 &:= \{(x,y) \mid  
y = \frac{z + \sqrt{z^2 -4}}{2} x \}.\notag
\end{align}
When $z < -2$, the asymptotes $l_1$ and $l_2$ lie entirely
in  quadrant II and IV, and
\[
\hcz \cap a_1 = \varnothing , \quad \hcz \cap a_2 \neq \varnothing
\]
When $z > 2$, the asymptotes $l_1$ and $l_2$ lie entirely
in  quadrant I and III, and
\[
\hcz \cap a_1 \neq \varnothing , \quad \hcz \cap a_2 = \varnothing
\]
\end{lemma}

\begin{figure}[htbp]
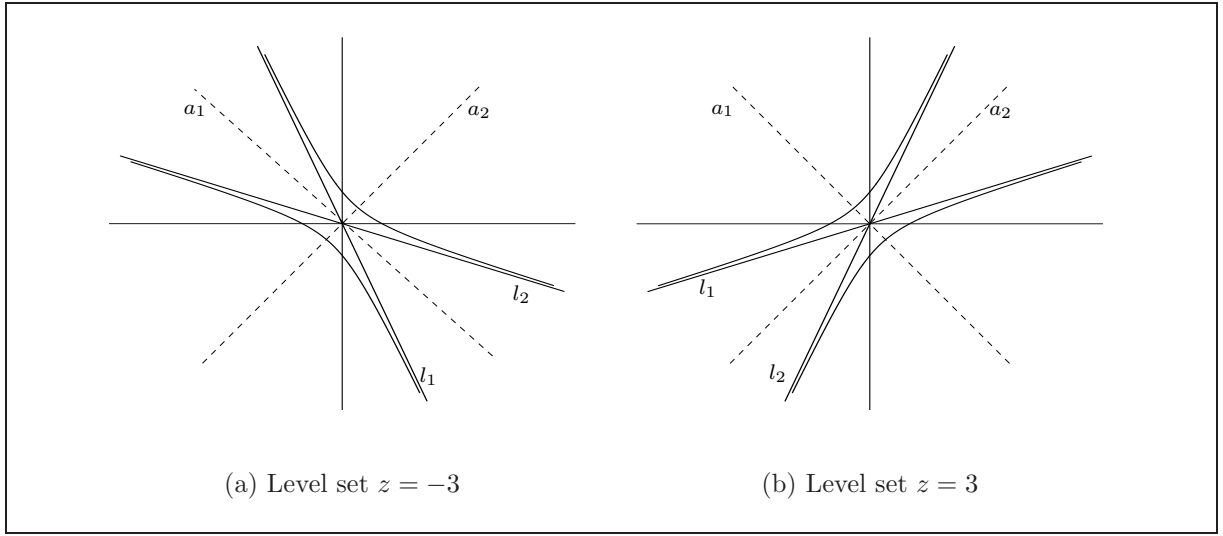

\centering
\subfigure[Level set $z=-3$]{\label{F:hypzneg}
\input{hypzneg.pstex_t}
}
\quad
\subfigure[Level set $z=3$]{\label{F:hypzpos}
\input{hypzpos.pstex_t}
}
\caption{Intersections of $\kappa\inv(c)$ with $z$-coordinate level sets}
\end{figure}
\begin{proof}
Assume $c<2$ and fix $z_0$ such that $|z_0|>2$. The
quadratic form
\[
S_{z_0}(x,y) = -x^2 -y^2 + z_0 xy
\]
has discriminant $D_{z_0} = 1 - z_0^2/4$, which is negative for $|z_0|>2$.
Therefore  $S_{z_0}$ is indefinite and its level sets must be
hyperbolae. The eigenvalues of $S_{z_0}$ are 
\[
s_1 = -\frac{z_0}{2} -1 \quad \text{and} \quad s_2 = \frac{z_0}{2} -1 
\]
with eigenvectors 
\[
\mathbf{e}_1 =\left(\begin{matrix} \frac{-1}{\sqrt{2}} \\ \frac{1}{\sqrt{2}}
\end{matrix}\right) \quad\quad\quad
\mathbf{e}_2 = \left(\begin{matrix} \frac{1}{\sqrt{2}} \\ \frac{1}{\sqrt{2}} 
\end{matrix}\right)
\]
Therefore with respect to the orthonormal basis 
$(\mbf{e}_1, \mbf{e}_2)$, the form $S_{z_0}$ has the following presentation
\[
S_{z_0}(\t{x}, \t{y}) =
-(\frac{z_0}{2} +1)\t{x}^2 + (\frac{z_0}{2} -1)\t{y}^2
\]
The principal axes $a_1$ and $a_2$ of $\hcz$ are, by definition,
collinear with $\mbf{e}_1$ and $\mbf{e}_2$ respectively, which implies
the first part of Lemma \ref{L:hypcz}
\par Since by assumption $z_0^2 -c -2 >0$, the equation of $\hcz$ 
\[
S_{z_0}(\t{x}, \t{y}) + (z_0^2 - c -2) =0 
\]
in terms of $(\mbf{e}_1, \mbf{e}_2)$, has the following canonical form
\begin{equation}\label{E:hypcan}
\left\{
\begin{aligned}
-\frac{\t{x}^2}{a_c^2(z_0)} + \frac{\t{y}^2}{b_c^2(z_0)} &=1, &
\text{if}\quad z_0<-2\\  
\frac{\t{x}^2}{a_c^2(z_0)} -\frac{\t{y}^2}{b_c^2(z_0)} &=1, &
\text{if}\quad z_0> 2 
\end{aligned} \right.
\end{equation}
where
\[
\begin{aligned}
a_c(z_0) &= \sqrt{2\left|\frac{(-z_0^2 + c +2)}{z_0 + 2}\right|} \\
b_c(z_0) &= \sqrt{2\left|\frac{(-z_0^2 + c +2)}{z_0 - 2}\right|}
\end{aligned}
\]
\par
Therefore in $\t{x}, \t{y}$-coordinates, the asymptotes of $h_c(z_0)$ are
\[
\t{l}_{1,2}: a_c(z_0)\t{y} \mp b_c(z_0)\t{x} = 0
\]
from which equations (\ref{E:equasmp}) follow after applying  the
appropriate coordinate transformation. Similarly, 
the principal axes of $h_c(z_0)$ are:
\[
\begin{array}{ll}
a_1: \t{y} = 0,  & a_2: \t{x} = 0
\end{array}
\]
implying \eqref{E:principalaxes}.
The assertions about the intersections 
$h_c(z_0)\cap a_1$ and $h_c(z_0)\cap a_2$ follow from \eqref{E:hypcan}.
\end{proof}
Next consider the restriction of $\tau$ to
$\hcz$. 
\begin{defin}
Let $l$ be a line in $\R^2$ with equation $l(x,y)=0$. Then $l$
partitions $\R^2$ into half-planes
\[
\begin{aligned}
H_l^+ &:= \{(x,y)\in\R^2\mid l(x,y)>0\} \\
H_l^- &:= \{(x,y)\in\R^2\mid l(x,y)<0\} \\
\end{aligned}
\]
the \emph{positive}, and the \emph{negative
half-plane associated with $l$} respectively. 
\end{defin}
\begin{lemma}\label{L:tauhyp}
Suppose $c<2$, $|z|>2$ and let $\hcz$ be as in Lemma
\ref{L:hypcz}. Let
\[
\mbf{a} := \left\{ 
\begin{array}{ll}
a_1, & \text{\rm ~if~}\, a_1\cap \hcz \neq \varnothing \\
a_2, & \text{\rm ~if~}\, a_2\cap \hcz \neq \varnothing 
\end{array}
\right.
\]
and $H_{\mbf{a}}^+$ and $H_{\mbf{a}}^-$ as above.
Then 
\begin{enumerate}
\item The restriction $\tau|_{\hcz}$ is invariant with respect to
reflections in $a_1$ and $a_2$.
\item  Let $\{P_1, P_2\} = \hcz \cap \mbf{a}$.
Then $\tau|_{\hcz}$ has global minima at $P_1, P_2$
\item $\tau|_{\hcz}$ decreases along $\hcz\cap H_{\mbf{a}}^-$ and
increases along $\hcz\cap H_{\mbf{a}}^+$.
\end{enumerate}
\end{lemma}
\begin{proof}
With respect to the standard basis in $\R^2$
reflection in $a_1$ and $a_2$ respectively are given by matrices  
\begin{equation*}
\bmatrix 
-1 &  0 \\
0 & -1 
\endbmatrix, \qquad 
\bmatrix 
0 &  1 \\
1 & 0 
\endbmatrix.
\end{equation*}
Each preserves $\zbar = -xy - z$ as well as $\tau =-z\zbar$.
\par 
Next assume $z < -2$. Parametrize \hcz\ as:
\begin{equation}\label{E:hypparzneg}
h_{c,z}(t)=
\left[
\begin{matrix}
x(t) \\  y(t)
\end{matrix}
\right] = 
\frac{\sqrt{2}}{2}
\left[
\begin{matrix}
-a_c(z) \sinh t \pm b_c(z)\cosh t \\
 a_c(z) \sinh t \pm b_c(z)\cosh t
\end{matrix}\right],
\quad t \in\R.
\end{equation}
For fixed $z$,  the restriction $\tau|_{\hcz}$  is a function of $t$ alone:
\[
\tau|_{\hcz}(t) = -z\zbar = x^2 + y^2 + c  + 2 = b_c^2(z)\cosh ^2t + 
a_c^2(z)\sinh ^2t + c + 2. 
\]
The sign of the derivative of $\tau$
\[
\frac{d\tau}{dt} = \sinh 2t(a_c^2(z) + b_c^2(z))
\]
depends only on the sign of $\sinh 2t$. Thus
\[
\frac{d\tau}{dt} \; \left\{
\begin{array}{lr}
< 0 & \text{if and only if}\; t<0 \\
= 0 & \text{if and only if}\; t=0 \\
>0  & \text{if and only if}\; t>0 
\end{array}
\right.
\]
Since
\[
\hcz\cap H_{\mbf{a}}^- = \{\hcz(t) \mid t < 0\}, \quad 
\hcz\cap H_{\mbf{a}}^+ = \{\hcz(t) \mid t > 0\} 
\]
and $\{P_1,P_2\} = \hcz(0)$, the last two parts of  Lemma~\ref{L:tauhyp}
follow. A similar argument applies when $z>2$.
\end{proof}
We complete the proof of Proposition~\ref{P:tautree}.
Assume $c<2$ and $|z_0|>2$. Let $u_0=(x_0, y_0, z_0)$ lie on \hczz. 
Since $Q_x$ interchanges the two (necessarily real)
roots of the quadratic polynomial
\[
\kappa(x,y_0,z_0) -c = -x^2 +(y_0z_0)x -y_0^2 +z_0^2 - 2 -c
\]
$Q_x$ acts as a deck transformation of the double covering
$\hczz\longrightarrow\R$.. 
Thus $l_x=\lyz\cap \lzz$ intersects \hczz\ precisely at
two points: $u_0$ and $Q_x(u_0)$. Similarly $l_y = \lxz\cap\lzz$
intersects \hczz\ at $u_0$ and $Q_y(u_0)$. 
\begin{claim} The points $u_0$ and $Q_x(u_0)$ lie on opposite branches of
\hczz. The points $u_0$ and $Q_y(u_0)$ lie on opposite branches of
\hczz.
\end{claim}
The line at infinity $\linf$ in
the projective completion $\rpt$ 
of the affine plane $L_{x-x_0}$ intersects 
$\hczz$ precisely at two points, the points of tangency
of \hczz with its asymptotes $l_1$ and $l_2$.
Denote these points by $\pinf_1$ and
$\pinf_2$. The point $O=l_1\cap l_2$ corresponds to the
origin $(0,0,z_0)$ in $L_{z-z_0}$. The coordinate 
axes correspond to lines $a_x$ and $a_y$ respectively, both of which
must intersect at $O$. Let $\pinf_x=\linf\cap a_x$ and
$\pinf_y=\linf\cap a_y$. Since in the affine plane  $\hczz$ intersects
both the $x$- and $y$-axis for any $|z_0|>2$, it intersects $a_x$ and
$a_y$ in $\rpt$. Therefore  $\pinf_x$ and $\pinf_y$ lie in the interior 
of the segment $\pinf_1\pinf_2$ 
(see Figure~\ref{F:projhyp}).  
$\linf$ partitions $\hczz$ into two disjoint arcs corresponding to the two
branches of $\hczz$ in the affine plane. 
A line in  $\R^2$ intersects both branches of
$\hczz$ if and only if the corresponding line in $\rpt$ intersects $\linf$
in the interior of the segment $\pinf_1\pinf_2$. 
But $u_0$ and $Q_x(u_0)$ span a line parallel to the $x$-axis, so 
its corresponding line in $\rpt$ intersects $a_x$ at $\pinf_x$. 
Similarly,  $u_0$ and $Q_y(u_0)$ span a line 
in $\rpt$ intersecting  $a_y$ at $\pinf_y$. 
The proof of the claim is complete..

Therefore the action of 
\[
\Lambda_{x,y} = \langle Q_x, Q_y \rangle \subset \Lambda
\]
on the plane $\lzz\subset\cv{1}{1}$ reduces to
a linear action of the infinite dihedral group on hyperbolae 
$\hczz$ (see Figure \ref{F:lxy-action}). 


\begin{figure}
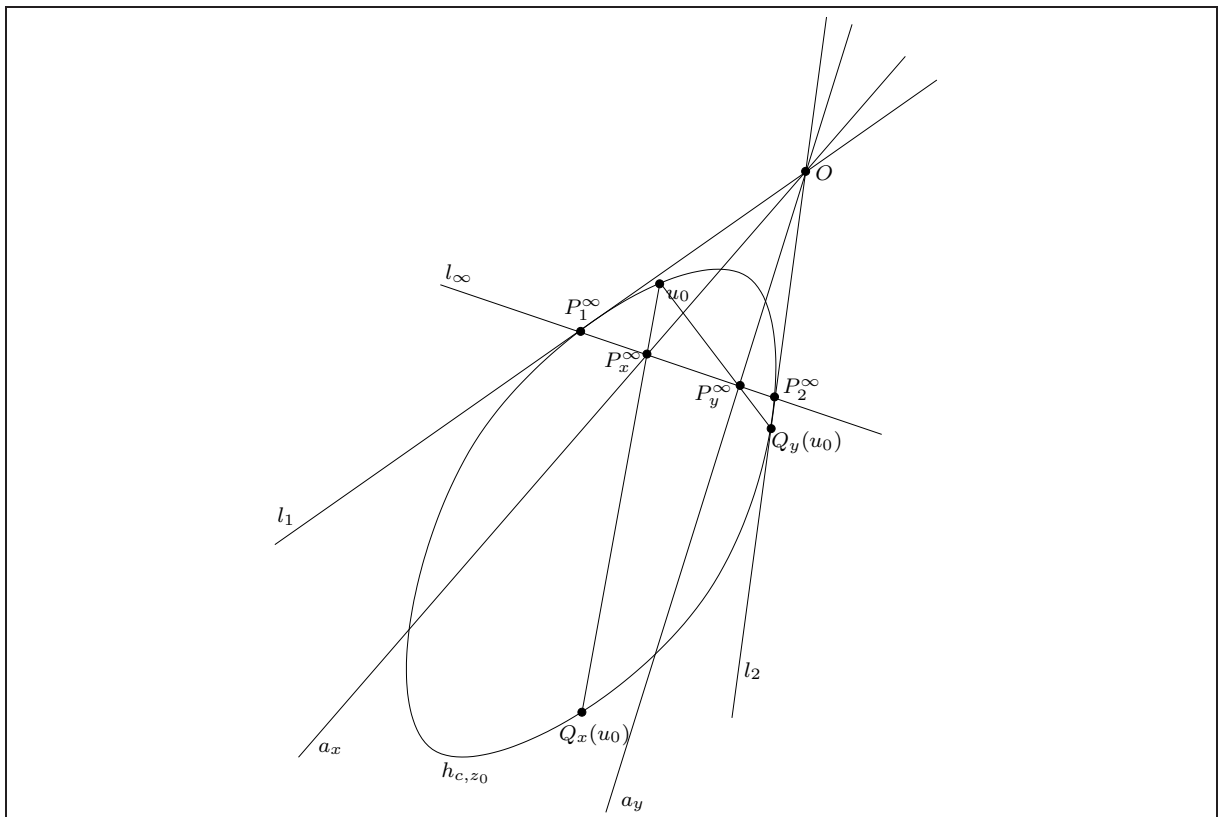

\centering
\input projhyp.pstex_t
\caption{Projective model of $\hczz$ and associated objects}
\label{F:projhyp}
\end{figure}

\begin{figure}[htbp]
\centering
\input{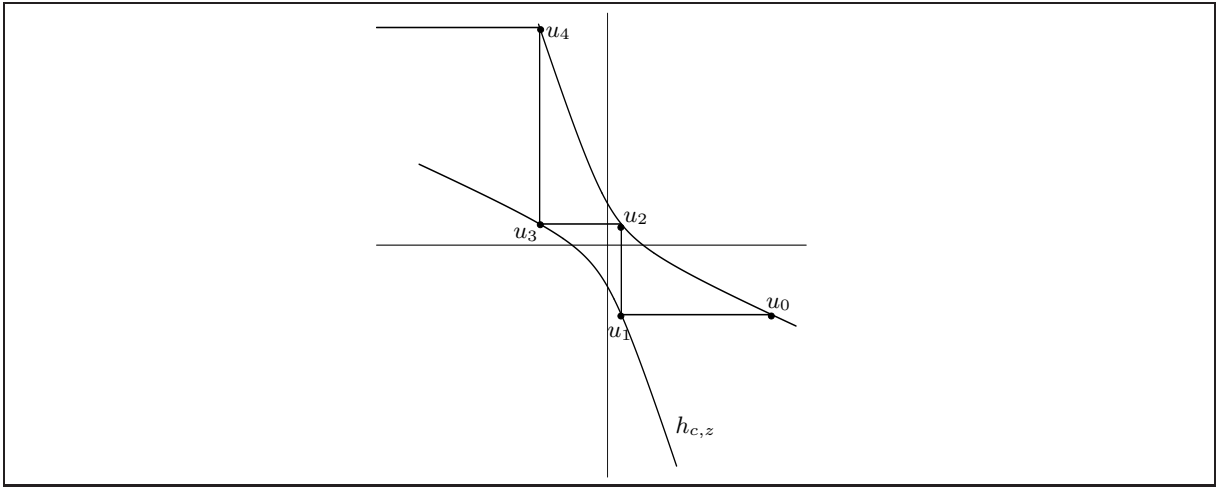}
\caption{Points on the $\Lambda_{x,y}$-orbit: 
$u_0$, $u_1 = Q_x(u_0)$, $u_2 = Q_y(u_1)$, etc}\label{F:lxy-action}
\end{figure}
\begin{defin}
We call the sets
\[
\begin{aligned}
\E{O}_x(u) &= \{u, Q_x(u), Q_y Q_x(u) \dots \}\\
\E{O}_y(u) &= \{u, Q_y(u), Q_x Q_y(u) \dots \}
\end{aligned}
\]
respectively the \emph{$x$-forward} and the \emph{$y$-forward
$\Lambda_{x,y}$-orbit} of $u$. 
\end{defin}
\pagebreak
Lemma~\ref{L:tauhyp} implies that  points on \hcz\
may be partitioned into four types 
\begin{enumerate}
\item type \tpp: points $u\in\hcz$ such that 
\[
\tau(Q_x(u)) - \tau(u)>0 , \quad \tau(Q_y(u)) - \tau(u) >0
\]
\item type \tpm: points $u\in\hcz$ such that
\[
\tau(Q_x(u)) - \tau(u) > 0 , \quad \tau(Q_y(u)) - \tau(u) <0
\]
\item type \tmp: points $u\in\hcz$ such that
\[
\tau(Q_x(u)) - \tau(u)<0 , \quad \tau(Q_y(u)) - \tau(u)>0
\]
\item type \tzz: points $u\in\hcz$ such that
\[
\begin{array}{cccc}
\text{either} & \tau(Q_x(u)) - \tau(u)=0, & \text{or} & 
\tau(Q_y(u)) - \tau(u)=0
\end{array}
\]
\end{enumerate}
Let $\E{T}$ denote the set $\{\text{\tpp, \tpm, \tmp, \tzz}\}$ of
point types. Define a function $\theta :\hcz \lrarw \E{T}$ mapping
each $u\in \hcz$ to the element of $\E{T}$ representing the type of
$u$. We call $\theta$ the \emph{type assignment function} on \hcz.
\par
Lemmas~\ref{L:hypcz},\ref{L:tauhyp} provide a classification for
the subsets of points of each type along~\hcz.  Exactly four 
points have type \tzz: the $x$- and the $y$-intercepts of
\hcz. Let $V_{l_i,x}$ (respectively $V_{l_i,y}$) denote the double
cone spanned by the asymptote $l_i$ and the $x$-axis (respectively $y$-axis)%
\footnote{There are, naturally, two double cones associated with a pair
of intersecting lines; we choose the one with the smaller cone
angle\label{F:cones}}. 
Assume $z<-2$. Then all points of type \tpp\ lie in 
quadrant I~and~III. Points of type \tpm\ and \tmp\ lie inside
$V_{l_1,y}$ and $V_{l_2,x}$ respectively. Assume $z>2$. Then all points of
type \tpp\ lie in quadrant II~and~IV. Points of type \tpm\ and \tmp\ lie inside
$V_{l_2,y}$ and $V_{l_1,x}$ respectively.
Clearly, reflection in the principal axes $a_1$ and $a_2$ interchages
the \tpm\ and \tmp\ types and leaves the \tpp\ type invariant.
A diagram representing points of different
types is shown on Figure~\ref{F:p-types-zneg} for the case $z< -2$ and on
Figure~\ref{F:p-types-zpos} for the case $z>2$. 
\begin{figure}[htbp]
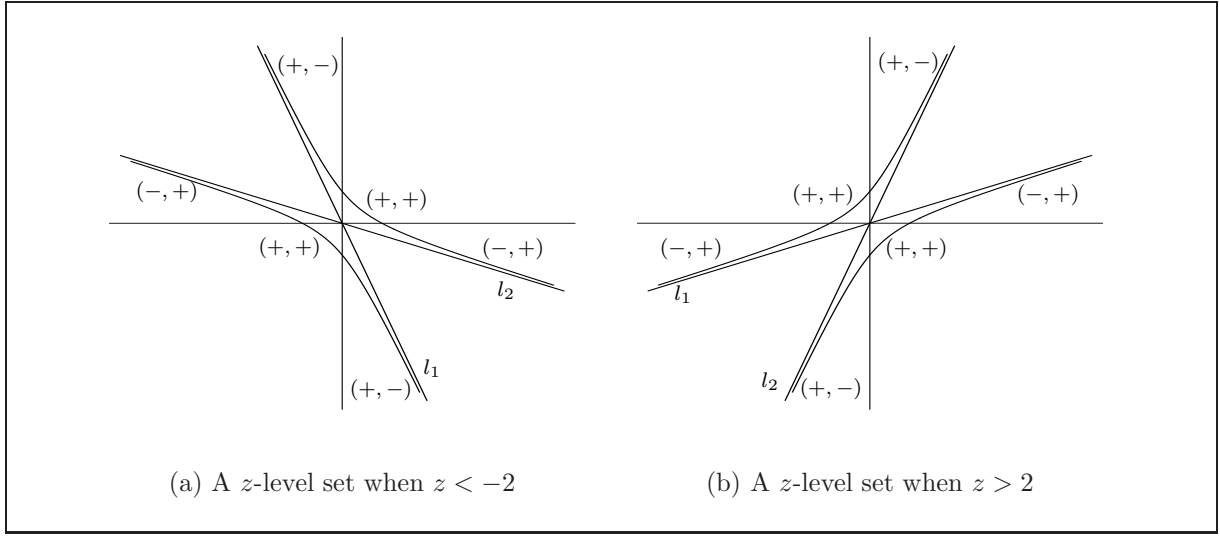

\centering
\subfigure[A $z$-level set when $z<-2$]{\label{F:p-types-zneg}
\input{point-types.pstex_t}
}
\quad
\subfigure[A $z$-level set when $z>2$]{\label{F:p-types-zpos}
\input{p-types-zpos.pstex_t}
}
\caption{Classification of Points With Respect to the Monotonicity of
$\tau$}
\end{figure}
\par
Lemma~\ref{L:omegaM} implies that \oM\ lies in the cylinder
\[
x^2 + y^2 < -c - 6, \; z \in \R
\]
Therefore, if for some $z<-2$ and $c< -14$, 
\[
L_z\cap\oM\cap\kappa\inv(c) \neq \emptyset,
\] 
then $L_z\cap\oM\cap\kappa\inv(c)$ consists of two arcs of 
\hcz\ inside the disk 
\[
D_c^{\Omega}: \quad x^2 + y^2 < -c-6.
\]

\begin{figure}[htbp]
\centering
\input{omegaDisk.pstex_t}
\caption{The set $L_z\cap\oM\cap\kappa\inv(c)$ inside $D_c^{\Omega}$}
\end{figure}

The $x$- and the $y$-intercepts of \hcz\  are:
\[
x_{1,2}^h = y_{1,2}^h = \pm\sqrt{z^2 - c -2}
\]
while the values of the $x$ and the $y$-intercepts of $\partial
D_c^{\Omega}$ are:
\[
x_{1,2}^D = y_{1,2}^D = \pm\sqrt{-c-6} \\
\]
Since $\sqrt{z^2 - c -2} > \sqrt{-c -6}$, the arcs $\hcz \cap
D_c^{\Omega}$ do not intersect the coordinate axes. Therefore all
points in $\hcz \cap D_c^{\Omega}$ have type \tpp, and
consequently 
\[
\tau(u) < \tau(Q_x(u)), \quad  \tau(u) < \tau(Q_y(u)) ,\quad
\forall u\in\hcz \cap D_c^{\Omega}
\]
Thus, when $u\in\oM\cap\kappa\inv(c)$, $\tau$ strictly increases
at the initial point  of each orbit $\E{O}_x(u)$ and $\E{O}_y(u)$.
(Compare Proposition~\ref{P:tautree}.)
Furthermore suppose that  $u\in\hcz$ satisfies
\[
\theta(Q_x(u)) \neq \text{\tpp} \quad \text{and} \quad 
\theta(Q_y(u)) \neq \text{\tpp}.
\]
Then 
\[
\begin{aligned}
\theta(u) &= \text{\tmp} \Rightarrow \theta(Q_x(u)) =
\theta(Q_y(u)) = \text{\tpm} \\
\theta(u) &= \text{\tpm} \Rightarrow \theta(Q_x(u)) =
\theta(Q_y(u))= \text{\tmp} \\ 
\end{aligned}
\]
Therefore, when $u\in\oM\cap\kappa\inv(c)$ the image of $\E{O}_x(u)$
under $\theta$ is the sequence 
\[
\text{\tpp,\tmp,\tpm,\tmp,}\dots.
\]
The image of $\E{O}_y(u)$ under $\theta$ is the sequence
\[
\text{\tpp,\tpm,\tmp,\tpm}\dots
\]
Consequently, $\tau$ strictly increases along both $\E{O}_x(u)$ and
$\E{O}_y(u)$.  
\par 
We return to the the proof of Proposition \ref{P:tautree}. 
Let $v$ be an arbitrary node of \blmd. Let 
\[
P_{u,v} = \{v_0=u, v_1, v_2,\dots, v_N=v\}
\]
be the (unique) path from $u$ to $v$ in \blmd. 
Partition $P_{u,v}$ into subsets 
\[
U_0 = \{v_{i_0}=v_0,\dots, v_{i_1}\}, U_1 = \{v_{i_1+1}, \dots, 
v_{i_2}\}, \dots,
U_m = \{v_{i_m+1}, \dots, v_{i_{m+1}} = v_N\}
\]
such that $v_{i_j+1} = Q_z(v_{i_j})$ and $U_j$ is a subset of either
 $\E{O}_x(v_{i_j+1})$, or  $\E{O}_y(v_{i_j+1})$, for each $j=0,\dots,m$. 
Since  either $U_0 \subset \E{O}_x(u)$, or  $U_0 \subset \E{O}_y(u)$,
$\tau$ is strictly increasing along $U_0$. Furthermore, since $\tau$
 is $Q_z$-invariant, 
\[
\tau(v_{i_j}) = \tau(v_{i_j+1})
\]
for each $j=0,  \dots, m$. 
Also, by assumption $z(v_{i_1}) = z(v_0) < -2$, and $\zbar(v_{i_1}) =
\zbar(v_0) >2$. Consequently
\[
z(v_{i_1+1}) = z(Q_z(v_{i_1})) = \zbar(v_{i_1}) > 2
\]
Recall that the hyperbola $h_c(z_{i+1})$  has the same principal axes 
as $h_c(z_i)$, as prescribed by Lemma~\ref{L:hypcz}
\[
\begin{aligned}
a_1:\; x - y &=0 \\
a_2:\; x + y &=0
\end{aligned}
\]
However the intersection properties of $a_j$ with $h_c(z_i)$ and
$h_c(z_{i+1})$ are different, namely:
\[
a_1 \cap h_c(z_i) \neq \varnothing, \quad a_2 \cap h_c(z_i) = \varnothing
\]
while
\[a_1 \cap  h_c(z_{i+1}) = \varnothing, \quad 
a_2 \cap h_c(z_{i+1}) \neq \varnothing
\]
Thus all points of type \tpp\ on $h_c(v_{i_1+1})$ lie 
in quadrant II and IV. Recall that $v_{i_1}$ is a point of type
either \tpm\ or \tmp\ and hence lies in quadrant $\mop{II}$ or $\mop{IV}$.
But then $v_{i_1+1}$, which is the image of $v_{i_1}$ under $Q_z$, must
also lie in quadrant $\mop{II}$ or $\mop{IV}$ and therefore it must be
of type \tpp. Hence $\tau$ must be strictly increasing along
$\E{O}_x(v_{i_1+1})$ and $\E{O}_y(v_{i_1+1})$, and therefore along
$U_1$. Similar argument applies to $U_j$, for $1 < j \leq m$. 
The proof of Proposition \ref{P:tautree} is complete.
\end{proof}

\subsection{Growth of the $\tau$ function}

So far we have shown that for each $u\in \oM$, the $\tau$ function is
non-decreasing along the path $P_{u,v}$ from $u$ to an arbitrary node
$v$ in \blmd. Moreover, $\tau$ is strictly increasing along $P_{u,v}$
except for a (possibly empty) subset 
of nodes at which $\tau$ is constant. Next, we estimate the variation
of $\tau$ 
among any pair of characters 
\[(w_0, w_1) \in \Delta_c :=\{(w_0, w_1) \in \kappa\inv(c) \times
\kappa\inv(c) \mid w_1 = \lambda w_0, \; \lambda\in \{Q_x, Q_y\}\}
\]
such that $|z(w_0)|>2$. 
\begin{prop}\label{P:taugrowth}
Suppose $c<2$. Let $w_0 = (x_0,y_0,z_0) \in \kappa\inv(c)$ be such that
$|z_0|>2$. Then
\begin{enumerate}
\item\label{I:taux} if $w_1 = (x_1, y_0, z_0) = Q_x(w_0)$
\begin{equation*} 
|\tau(w_1) - \tau(w_0)| \geq |y_0z_0|
\max \left(|y_0|\sqrt{z_0^2 - 4},  2\sqrt{z_0^2 - c -2}\right) 
\end{equation*}
\item\label{I:tauy} if $w_1 = (x_0, y_1, z_0) = Q_y(w_0)$
\begin{equation*}
|\tau(w_1) - \tau(w_0)| \geq |x_0z_0|
\max \left(|x_0|\sqrt{z_0^2 - 4}, 2\sqrt{z_0^2 - c -2}\right)
\end{equation*}
\end{enumerate}
\end{prop}
\begin{proof}
Since $|z_0| > 2$ , the $z_0$-level set $L_{z-z_0}\cap \kappa\inv(c)$ is
the hyperbola \hczz. We prove case~(\ref{I:tauy}) first.
Let $A_{l_1}$ and $A_{l_2}$ be the  intersection points of the line
$l_x: x=x_0$ with $l_1$ and $l_2$ respectively. 
Clearly
\[
|y_1 - y_0| > \mop{d}(A_{l_1}, A_{l_2}) = 
\left|\frac{z_0 + \sqrt{z_0^2 -4}}{2}x_0 -
\frac{z_0 - \sqrt{z_0^2 -4}}{2}x_0\right| = |x_0|\sqrt{z_0^2 - 4}
\]
where $\mop{d}(\cdot\;,\cdot)$ denotes the (Euclidean) distance function in
$\R^2$. On the other hand, since the minimum vertical distance between
pairs of points on \hczz\ is achieved at $x_0=0$
\[
|y_1 - y_0| > \min_{w\in\hczz}|y(w) - y(Q_y(w))| = 2\sqrt{z_0^2 - c - 2}
\]
Therefore
\begin{equation}\label{E:yincr}
|y_1 - y_0| \geq \max \left(|x_0|\sqrt{z_0^2 - 4}, 
2\sqrt{z_0^2 - c -2} \right)
\end{equation}
By the defintion of $\tau$
\[
\begin{split}
|\tau(w_1) - \tau(w_0)| = &|z_0||(\zbar_1 - \zbar_0)| \\
= &|z_0||y_0||x_1 - x_0|
\end{split}
\]
which together with inequality~(\ref{E:yincr}) implies part
(\ref{I:tauy}) of Proposition \ref{P:taugrowth}. 
Similar argument applies when $w_1 = Q_x(u)$, in which case
\begin{equation}\label{E:xincr}
|x_1 - x_0| \geq \max \left(|y_0|\sqrt{z_0^2 - 4}, 
2\sqrt{z_0^2 - c -2}\right)
\end{equation}
and 
\[
\begin{split}
|\tau(w_1) - \tau(w_0)| = &|z_0||(\zbar_1 - \zbar_0)| \\
= &|z_0||x_0||y_1 - y_0|
\end{split}
\]
\end{proof}
The assumption that $|z|>2$ and $c<2$ eliminates $z$ from
the estimates for the increment of $\tau$. 

\begin{cor}\label{C:taubound}
Let $c<2$ be fixed, and let
\[
(w_0, w_1)\in \Delta_c
\]
such that $|z(w_0)|>2$, and $|z(w_1)|>2$. Then the variation of $\tau$ 
is bounded from below by a quantity depending on $x_0$ or $y_0$
alone. In particular 
\begin{equation}
|\tau(Q_x(w_0)) - \tau(w_0)| \geq 16|y_0|\sqrt{2-c}
\end{equation}
and
\begin{equation}
|\tau(Q_y(w_1)) - \tau(w_0)| \geq 16|x_0|\sqrt{2-c}
\end{equation}
\end{cor}

\subsection{$\tau$-reduction for Fricke-space characters}

We have seen so far that the $\tau$-function is non-decreasing 
along the levels of 
\blmd\ for any 
\begin{equation*}
u = (x, y, z) \in\oM\cap \kappa\inv(c). 
\end{equation*}
Also for every $w\in \blmd\cap \Delta_c$, 
and fixed $c$, the increment is bounded from below by a 
quantity depending on $|x|$ or $|y|$ alone. Since
$z<-2$ and $\zbar<-2$,
\begin{align*}
z\zbar & = -x^2 - y^2 - 2 -c > 4 \\
z + \zbar & = -xy < -4.
\end{align*}
The first inequality implies that $|x| < \sqrt{-c-6}$ and $|y| <
\sqrt{-c-6}$, which together with the second inequality yield
\[
\begin{aligned}
|x| & > \frac{4}{\sqrt{-c-6}} \\
|y| & > \frac{4}{\sqrt{-c-6}}.
\end{aligned}
\]
Thus if $\kappa\inv(c)\cap \oM\neq\varnothing$ (or equivalently, if $c<
-14$), there is a constant $T_c$ such that 
\begin{equation}
|\tau(w_1) - \tau(w_0)| \geq T_c
\end{equation}
for every $u\in\oM$ and a pair $w_0, w_1$ of successive nodes in
$\blmd\cap\Delta_c$. The uniform lower bound for the increment of
$\tau$ along the depth levels of \blmd\ guarantees the following. The
path $P_{v,u}$ from an arbitrary character $v$ in 
$\Lambda\cdot\oM$ to the unique
character  $u\in\oM$ such that $v = \blmd$, can be recovered in
finite number of steps via reduction of $\tau$ by a definite amount at
each step. More precisely,
\begin{lemma}\label{L:taured}
Let  $c<-14$ and let $v\in\kappa\inv(c)$. Assume $v \in \blmd$ for
some $u\in\oM$. Then 
\[
P_{v,u} = \{v_0=v, v_1, \dots, v_n =u\}
\]
where $\{v_i\}_{i=0}^n$ is the (unique) sequence such 
that for each $0\leq i \leq n-1$
\[
v_{i+1} = \left\{ 
\begin{array}{ll}
Q_x(v_i) & \text{if} \quad\tau(Q_x(v_i)) < \tau(v_i) \\
Q_y(v_i) & \text{if} \quad \tau(Q_y(v_i)) < \tau(v_i) \\
Q_z(v_i) & \text{otherwise} 
\end{array}
\right.
\]
Moreover $n$ depends only on $c$ and $v$. 
\end{lemma}
\pagebreak

\subsection{$\tau$-reduction for arbitrary characters}
This type of reduction process extends to arbitrary characters
\[
v \in \E{H}_c :=\kappa\inv(c) \cap
\big(\cv{1}{1} - \oM \cup \ec \cup Q_z(\ec)\big)
\]

By definition, $v\in \E{H}_c$ implies that  $|z(v)|>2$ and
$|\zbar(v)|>2$. 
For such $v$, the proof of Proposition \ref{P:tautree}
produces a sequence 
\[
\mcv: \quad v_0=v, v_1, \dots, v_n, \dots
\]
such that for each $i=1,2,\dots$: 
\begin{itemize}
\item $v_{i+1} = \lambda_i v_i$, where $\lambda_i\in\{Q_x, Q_y, Q_z\}$
\item $\lambda_i \neq \lambda_{i+1}$
\item $\tau(v_{i+1}) \leq \tau(v_i)$
\end{itemize}
\begin{defin}
We call such a sequence \emph{$\tau$-minimizing}. A $\tau$-minimizing
sequence is said to \emph{terminate}, if there exists $n$, such that
one of the following occurs
\begin{itemize}  
\item $v_n\in\ec$, or
\item $v_n \in \oM$, or
\item $v_n\in\scz:=\{(x,y,z)\in \kappa\inv(c)\mid x=0,\; \text{or}\; y=0\}$
\end{itemize}
The element $v_{n_t}$ indexed by the smallest such n, will be called
\emph{a terminator}. 
\end{defin}
The sequence $\mcv$ can be constructed as follows. Let $v_i=(x_i,
y_i, z_i)$, and assume $|z_i|>2$ and $|\zbar_i|>2$. Then $v_i$ lies 
on the hyperbola $h_c(z_i)$. There are four possibilities based on the
type of $v_i$
\begin{enumerate}
\item $v_i$ is of type \tmp; then $v_{i+1} := Q_x(v_i)$\label{I:xmin}  
\item $v_i$ is of type \tpm; then $v_{i+1} := Q_y(v_i)$\label{I:ymin} 
\item $v_i$ is of type \tpp; then $v_{i+1} := Q_z(v_i)$
\item $v_i$ is of type \tzz; then $v_{i}\in\scz$ and is therfore a terminator.
\end{enumerate}
\begin{prop}\label{P:tauterm}
After a finite number of steps, every $\tau$-minimizing
sequence terminates.
\begin{proof}
Clearly, this is true if $v\in \Omega^M$ (Lemma~\ref{L:taured}). 
Thus assume that $v\notin \Omega^M$.
\par
Recall that at each $v_i$ where $\tau$ is strictly decreasing,
Proposition \ref{P:taugrowth} provides a lower bound on 
$|\tau(v_i) - \tau(v_{i+1})|$ that depends on $x_i$ or $y_i$ alone. 
As long as the $x_i$'s and the $y_i$'s do not accumulate at $0$, the
decrement of $\tau$ by a definite amount at each step guarantees that
after finitely many steps there will be an element $v_n$ such that
$\tau(v_n)<4$ and hence $v_n$ will lie in one of $\scz$, \ec, or
$Q_z(\ec)$. Thus it suffices to show that if $v \notin \Lambda\cdot\scz$, then 
the $x$ and the $y$ coordinates of characters on the $\Lambda$-orbit
of $v$ but 
outside of \ec\ are bounded away from $0$ by some positive constant
depending on $v$.  
\begin{lemma}\label{L:xyaccum}
For every $v\notin\Omega^M\cup\ec\cup\Lambda\cdot\scz$ there exists positive
constants $\eps_x(v)$ and $\eps_y(v)$ such that
\[
\begin{aligned}
\inf_{\lambda\in\Lambda} \{x(\lambda v)\mid \lambda v \notin \ec\} =
\eps_x(v) > 0\\ 
\inf_{\lambda\in\Lambda} \{y(\lambda v)\mid \lambda v \notin \ec\} =
\eps_y(v) > 0\\ 
\end{aligned}
\]
\begin{proof}
We prove this by contradiction. Suppose no such constants exist.
If $|z|>2$ and $|\zbar|>2$, but $\hcz\cap \oM =\varnothing$, 
the set of points on \hcz\ of type \tpp\ partitions
into three disjoint subsets
\begin{align*}
\bar{Z}_x &= \{u \in \hcz \mid  Q_z(u) \in (-,+)\} \\
\bar{Z}_y &= \{u \in \hcz \mid  Q_z(u) \in (+,-)\} \\
\bar{Z}_e &= \{u \in \hcz \mid   -2 < z(Q_z(u)) < 2\} 
\end{align*}
\begin{defin}
Points in $\bar{Z}_x$ (respectively $\bar{Z}_y$, $\bar{Z}_e$) will be called
\emph{points of type} $\zbxneg$ (respectively $\zbyneg$,~$\zbell$). 
Notice that by definition $\bar{Z}_e \subset \ec$.
\end{defin}
Fix $c<2$ and consider the projection of $\kappa\inv(c)$ onto the
$xy$-plane
\[
\Pi: \kappa\inv(c)\ni (x,y,z) \longmapsto (x,y)\in \R^2.
\]
The level set $L_{z+2}\cap\kappa\inv(c)$ projects onto the pair
of lines
\[
h_{-2}^{\pm}: x + y \pm \sqrt{-c+2} =0
\]
and the set $L_{z-2}\cap\kappa\inv(c)$ projects onto the pair of lines
\[
h_{+2}^{\pm}: x - y \pm \sqrt{-c+2} =0.
\]
The pair $h_{-2}^{\pm}$ partitions the $xy$-plane into two regions:
\[
H_{-2}^{-} = \{(x,y)\in \R^2\mid (x+y +\sqrt{-c+2})(x+y
-\sqrt{-c+2})<0 \}
\]
and
\[
H_{-2}^{+} = \{(x,y)\in \R^2\mid (x+y +\sqrt{-c+2})(x+y
-\sqrt{-c+2})>0 \}
\]
\begin{figure}
\centering
\input{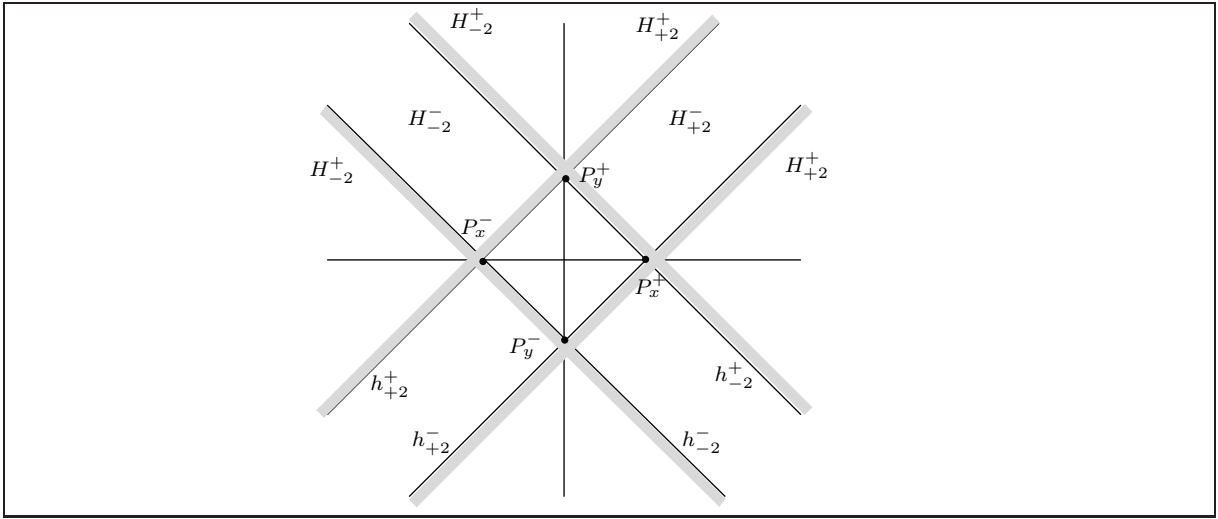} 
\caption{Lines, cones, and half-planes}\label{F:zbarcones}
\end{figure}
The first one, $H_{-2}^{-}$ is the ``strip'' bounded by $h_{-2}^{\pm}$, 
and the second one,  $H_{-2}^{+}$ is the complement of $H_{-2}^-$ in
$\R^2$ (see Figure \ref{F:zbarcones}).
Similarly, the pair $h_{+2}^{\pm}$ partitions the plane into regions
\[
H_{+2}^{-} = \{(x,y)\in \R^2\mid (x-y +\sqrt{-c+2})(x-y
-\sqrt{-c+2})<0 \}
\]
 and
\[
H_{+2}^{+} = \{(x,y)\in \R^2\mid (x-y +\sqrt{-c+2})(x-y
-\sqrt{-c+2})>0 \}
\]
When $z_0<-2$ (respectively $z_0>2$), the level sets $\lzz\cap
\kappa\inv(c)=\hczz$ 
project to a family of hyperbolae contained in $H_{-2}^+$
(respectively $H_{+2}^+$).  The regions $H_{-2}^-$ and $H_{+2}^-$ contain
a family of ellipses, which are the projections of the level sets 
$\lzz\cap \kappa\inv(c)$ for $-2< z_0 < 2$. Hence the set of
$\zbell$ points projects onto the union 
\[
H_{-2}^- \cup H_{+2}^-
\]
If $|z_0|>2$, the set of \zbxneg\ and \zbyneg points on \hczz\  projects onto 
\[
H_{-2}^+\cap H_{+2}^+
\]
which is a union of four cones
with vertices $P_x^{\pm} = (\pm\sqrt{-c+2}, 0)$ and 
$P_y^{\pm} = (0, \pm\sqrt{-c+2})$ (see Figure \ref{F:zbarcones} and \ref{F:zbarneg-pos}). 
\begin{figure}
\centering
\subfigure[for $z<-2$]{
\input{zbarneg.pstex_t}
}
\quad
\subfigure[for $z>2$]{
\input{zbarpos.pstex_t}
}
\caption{Points of type $\zbxneg$, $\zbyneg$, and $\zbell$}\label{F:zbarneg-pos}
\end{figure}
In particular, the set of \zbxneg\ points projects onto the
cones with vertices $P_x^{\pm}$. 
The set of \zbyneg\  points projects onto the cones with
vertices $P_y^{\pm}$ (see Figure~\ref{F:families}). 
\begin{figure}[hbtp]
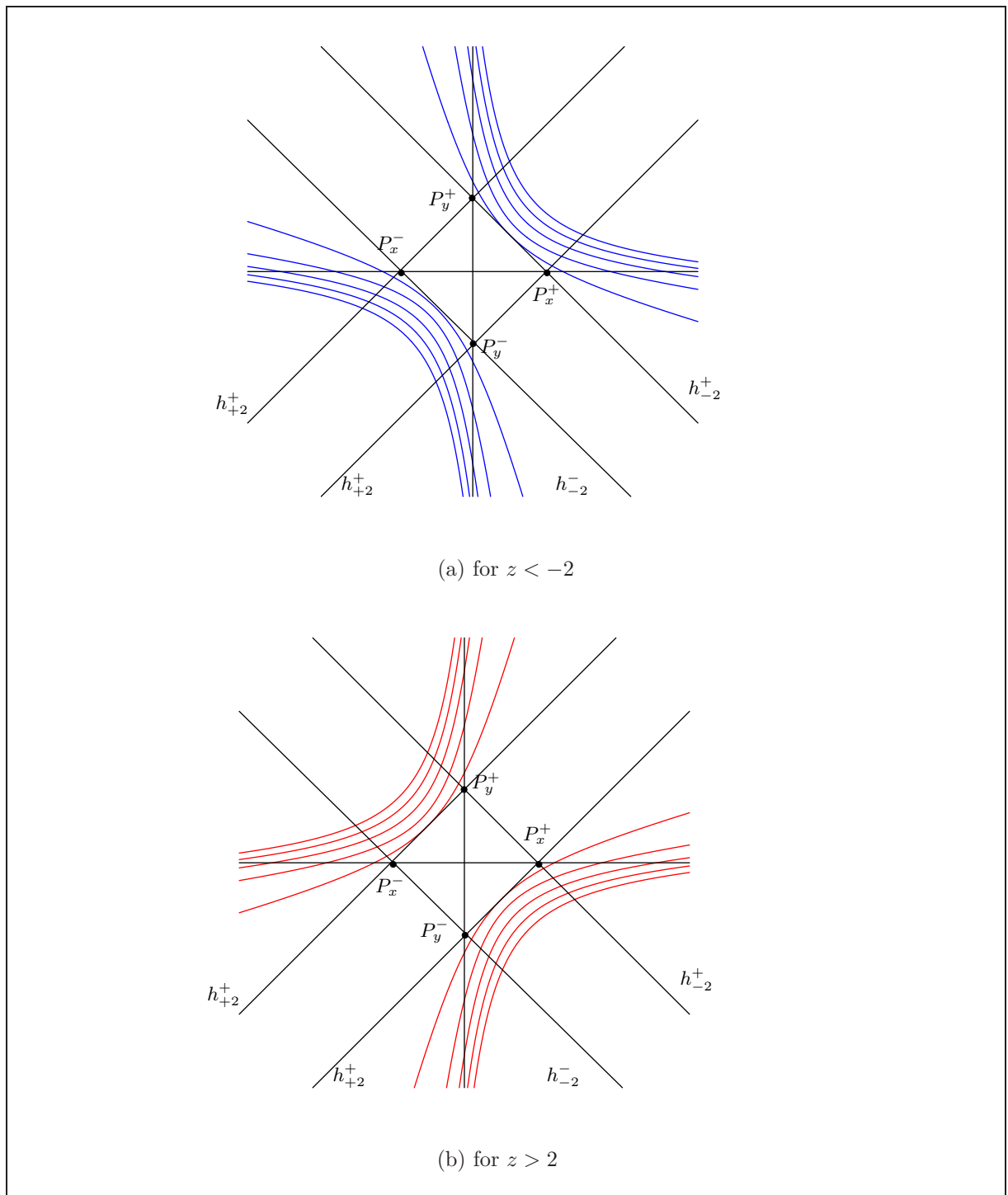

\centering
\subfigure[for $z<-2$]{
\input{hypfamneg.pstex_t}
}
\quad
\subfigure[for $z>2$]{
\input{hypfampos.pstex_t}
}
\caption{Families of hyperbolae as projections of $\kappa\inv(c)\cap
L_z$}\label{F:families} 
\end{figure}
Observe that the strip
\[
F_x = \R\times (-\sqrt{-c+2}, \sqrt{-c+2})
\]
does not contain any projections of \zbyneg\  points. 
Similarly, the strip
\[
F_y = (-\sqrt{-c+2}, \sqrt{-c+2})\times \R
\]
does not contain any projections of \zbxneg\  points. 
Consequently, if $v_{n_1}\notin\ec$ of the $\tau$-minimizing sequence
$\{v_i\}$ enters the slab
\[
F_x \times \R
\]
then it must remain on the same $y$-level set for all $i\geq n_1$.
Similarly, if an element $v_{n_2}\notin\ec$ enters the slab
\[
F_y \times \R
\]
then it must remain on the same $x$-level set for all $i\geq n_2$.
Since we had assumed that either the $x$ or the $y$-coordinates of
characters in $\{v_i\}$ accumulate at $0$, there exists
$v_n\notin\ec\cup\scz$, $n>0$, such that either
\[
v_n\in F_x \times \R
\]
or
\[
v_n\in F_y \times \R
\]
Suppose $v_n\in F_y \times \R$. This means that $x(v_i)$ must be constant
and $|y(v_i)|> \sqrt{-c+2}$ for all $i>n$, a contradiction.
Similarly, $v_n\in F_x \times \R$ implies a contradiction. This
concludes the proof of Lemma \ref{L:xyaccum} and hence of Proposition~\ref{P:tauterm}.
\end{proof}
\end{lemma}
\renewcommand{\qedsymbol}{}
The proof of Lemma \ref{L:xyaccum} suggests that if a 
$\tau$-minimizing sequence enters one of the regions $F_x\times\R$, or
$F_y\times\R$, its behavior can be determined explicitly. 
More precisely, let $\{v_i\}$, be a $\tau$-minimizing sequence. 
\begin{defin}
A \emph{terminal $x$-plane} (respectively $y$-plane) is a level set $L_{x-x_0}$  
(respectively $L_{y-y_0}$) such that if $v_n \in L_{x-x_0}$
(respectively $v_n \in L_{y-y_0}$) for some $n$, then $v_i \in
L_{x-x_0}$ (respecitvely $v_i \in L_{y-y_0}$), for all $n\leq i \leq
n_t$, where $n_t$ is the index of the terminating element of $\{v_i\}$. 
\end{defin}
Thus all elements of a $\tau$-minimizing sequence that belong to
\[
F_x\times\R\cup F_y\times\R
\]
must lie on a terminal plane. 
Figure \ref{F:termplane} shows several elements of the minimizing sequence of 
the character $u=(-0.2, 12, -10)$  on its terminal
$x$-plane. Successive points are joined by line segments to facilitate
visualization. 
\begin{figure}[htbp]
\centering
\input{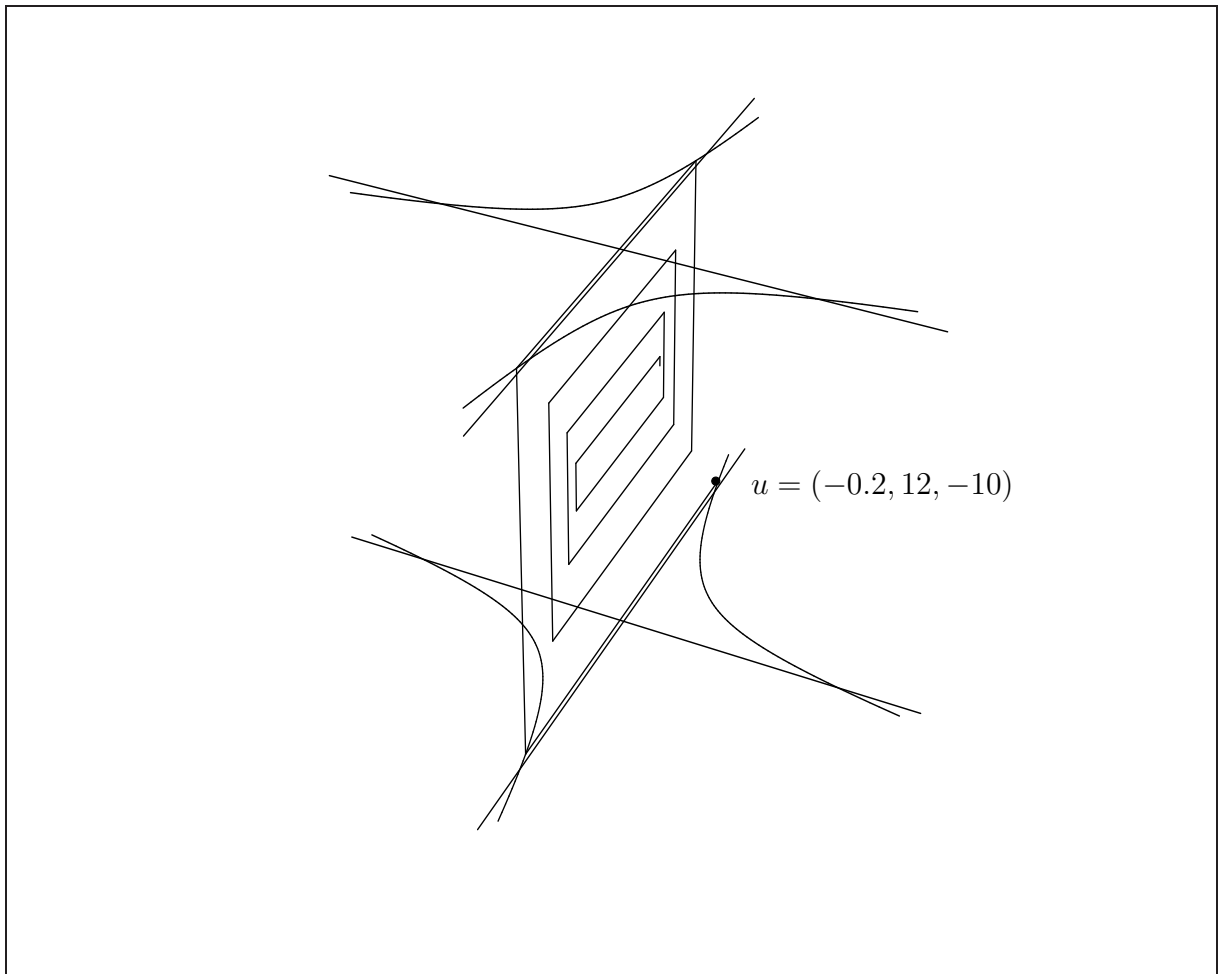}
\caption{The terminal plane of the character $u=(-0.2, 12, -10)$}\label{F:termplane}
\end{figure}
\newpage
\subsection{The $\tau$-Reduction Algorithm}
We use the results obtained so far to construct an algorithm, which
implements the method of $\tau$-reduction for distinguishing among
characters on $\kappa\inv(c)$ when $c<-14$.
\begin{tabbing}
\hspace*{.25in}\=\hspace{3ex}\=\hspace{3ex}\=\hspace{3ex}\kill
\> \textbf{Algorithm} \textsc{$\tau$-reduction}$(u)$ \\
\> \textit{Input.} A  character $u \in \kappa\inv(c)$. \\
\> \textit{Output.} A character in $\oM$, $\ec$, or $\scz$, which is 
$\Gamma$-equivalent to $u$. \\
\> $u\leftarrow u_0$; \\
\> \textbf{while} $|\zbar| > 2$ \textbf{do} \\
\> \> \textbf{if} x(u) = 0 \textbf{or} y(u) = 0 \textbf{then} $u \in
\scz$ \textbf{return}\ $u$ \\
\> \> \textbf{if} $z(u) < -2$ \textbf{and} $\zbar(u)< -2$ \textbf{then} $u \in
\oM$ \textbf{return}\ $u$ \\
\> \> \textbf{if} $z(u) > 2$ \textbf{and} $\zbar(u) > 2$ \textbf{then} 
$\sigma_{xz}(u) \in \oM$ \textbf{return} $\sigma_{xz}(u)$ \\
\> \> \textbf{if} $\tau(Q_x(u)) < \tau(u)$ \textbf{then} $u
\leftarrow Q_x(u)$ \\
\> \> \textbf{elseif} $\tau(Q_y(u)) < \tau(u)$ \textbf{then} 
$u \leftarrow Q_y(u)$ \\
\> \> \textbf{else} $u \leftarrow Q_z(u)$ \\
\> \textbf{end do} \\
\> $u\in \ec$; \textbf{return}\ $u$.
\end{tabbing}
\end{proof}
\end{prop}
\section{The Action of the Modular Group on Characters}
In Chapter~\ref{S:fricke} we proved that the action of $\Gamma$ on 
\[
\Omega^M = \coprod_{\gamma\in\Gamma}\gamma\oM
\]
is wandering. In this chapter we show that the action of $\Gamma$ on
the complement of $\Omega^M$ in $\cv{1}{1}$ is ergodic in the
following sense. Recall that the action of $\Gamma$ induces a measurable
equivalence relation ``$\sim$'', which is ergodic if and only if every
function that is constant on equivalence classes is constant almost
everywhere. In that context, if every point in a subspace $X$ is
$\Gamma$-equivalent to a point in a subspace $Y$, then ergodicity on $X$
is equivalent to ergodicity on $Y$ (regardless of whether $Y$ is
invariant or not; compare Goldman~\cite{Gold:pTorus}). In the previous
section we have proved
that when $c<2$ every point on
$(\cv{1}{1}-\Omega^M\cup\ec)\cap\kappa\inv(c)$
is $\Gamma$-equivalent to a point in \ec. Thus ergodicity
on $\cv{1}{1}$ reduces to ergodicity 
of ``$\sim$''  on \ec, or equivalently on its ``$Q_z$-dual''
\[
\ecb = Q_z(\ec) = \{(x,y,z)\in\kappa\inv(c)\mid -2 < z < 2\}
\]
Let $u=(x,y,z)\in \ecb$. Then
$e_c(z) = L_z\cap\kappa\inv(c)$ is an ellipse upon which $Q_xQ_y \in \Lambda_{x,y}$ acts
by a rotation of angle
$$
\alpha = 2\cos^{-1}\frac{z}2.
$$
(Compare Goldman~\cite{Erg}.) 
For almost every $z$ (namely when $\alpha/(2\pi)$ is irrational), 
this transformation generates an ergodic action
on $e_c(z)$. Thus a function
\[
f:\; \ecb \lrarw \R
\]
that is $\Lambda$-invariant, would be constant almost everywhere
on each $L_z\cap\kappa\inv(c)$ and would therefore depend
almost everywhere on $z$ alone. Thus $f(x,y,z) = g(z)$ almost
everywhere for some function 
\[
g:\; [-2,2] \lrarw \R
\]
It now suffices to eliminate the dependence of $g$ on $z$. To this
end, we parametrize $e_c(z)$ as follows (see Lemma \ref{L:hypcz})
\[
e_c(z):\; \left\{
\begin{aligned}
x &= \frac{\sqrt{2}}{2}(-A\cos\theta + B\sin\theta)\\
y &= \frac{\sqrt{2}}{2}(\phantom{-}A\cos\theta + B\sin\theta)
\end{aligned}
\right.
\]
where 
\[
A = \sqrt{\frac{z^2 -c -2}{2+z}}, \quad B= \sqrt{\frac{z^2 -c -2}{2-z}}
\]
Consequently, for a fixed $z\in (-2, 2)$ the restriction of $\zbar$ to
$e_c(z)$ depends on $\theta$ alone
\[
\zbar_{c,z}(\theta) = \zbar|_{e_c(z)} = \frac{z^2 - c- 2}{2+z}\cos^2\theta - 
\frac{z^2 -c -2}{2-z}\sin^2\theta - z
\]
whose extrema are attained at $\theta=k\pi/2$, for $k\in\Z$. The
critical values of  $\zbar_{c,z}(\theta)$ are 
\begin{equation}
\begin{aligned}
\zbar_{\text{odd}} &= -\frac{z^2 - c -2}{2-z} -z &= 
\frac{-2z +c +2}{2-z} &= \phantom{-}2 + \frac{c-2}{2-z}\\
\zbar_{\mop{even}}&= \phantom{-}\frac{z^2 - c -2}{2+z} -z &= 
\frac{-2z -c -2}{2+z} &= -2 - \frac{c-2}{2+z} \\
\end{aligned}
\end{equation}
respectively, depending on the parity of $k$. We classify the extrema
of $\zbar_{c,z}$ by analysing the second derivative
\[
\zbar_{c,z}''(\theta) = -\frac{8(z^2 -c -2)}{4-z^2}\cos2\theta 
\]
Clearly
\[
\zbar_{c,v}''\left(\frac{k\pi}{2}\right) \quad\left\{
\begin{array}{ll}
>0, &  k\; \text{odd} \\
<0, &  k\; \text{even}
\end{array}
\right.
\]
provided that 
\begin{equation}\label{E:zsq}
z^2 - c -2 >0
\end{equation}
Observe that when $c<-2$, the latter is 
satisfied trivially for any $z\in \R$. On the other hand, when 
$-2< c< 2$, inequality (\ref{E:zsq}) 
is satisfied \emph{a fortiori} for all points in
\[
\ecb = \kappa\inv(c) \cap \R^2\times(-2, 2)
\]
since for any such point $u=(x,y,z)$, the quadratic form 
\[
S_z(x,y)=-x^2 -y^2 +xyz
\]
is negative definite and therefore the set
\[
\kappa\inv(c)\cap \R^2\times(-\sqrt{c+2}, \sqrt{c+2})
\]
is empty. Consequently, the critical value $\zbar_{\text{odd}}$ is a minimum,
and the critical value $\zbar_{\text{even}}$ is a maximum of
$\zbar_{c,z}$ for each $z \in (-2, 2)$, such that $L_z \cap
\kappa\inv(c)$ is non-empty. Clearly, these extrema are global. 
Thus for any $z\in (-2 ,2)$
\begin{equation}
\begin{aligned}
\zbar_{min}(c,z) &= \min_{\theta\in\R} \zbar_{c,z}(\theta) &=  
\phantom{-}2 + \frac{c-2}{2-z} \\ 
\zbar_{max}(c,z) &= \max_{\theta\in\R} \zbar_{c,z}(\theta) &= 
-2 - \frac{c-2}{2+z}  
\end{aligned}
\end{equation}
\par
When $c<-14$ 
\begin{equation} \label{E:zminineq}
(-2, 2) \subsetneq (\zbar_{min}(c,z), \zbar_{max}(c,z))
\end{equation}
for any $z\in(-2,2)$. 
The cyclic group generated by an irrational rotation of an 
ellipse acts ergodically and therefore,  for almost $z\in (-2, 2)$ and
every $u\in e_c(z)$, the set  
\[
\{\zbar(\lambda u)\mid \lambda\in\Lambda_{x,y}\} 
\]
must be dense in $[\zbar_{min}(c,z), \zbar_{max}(c,z)]$.
Therefore the \mbox{$Q_z$-invariance} of $f$ implies that $f$
must be constant also with respect to $z$, for almost every $z\in (-2,2)$. 
\par
\begin{figure}[htbp]
\centering
\input{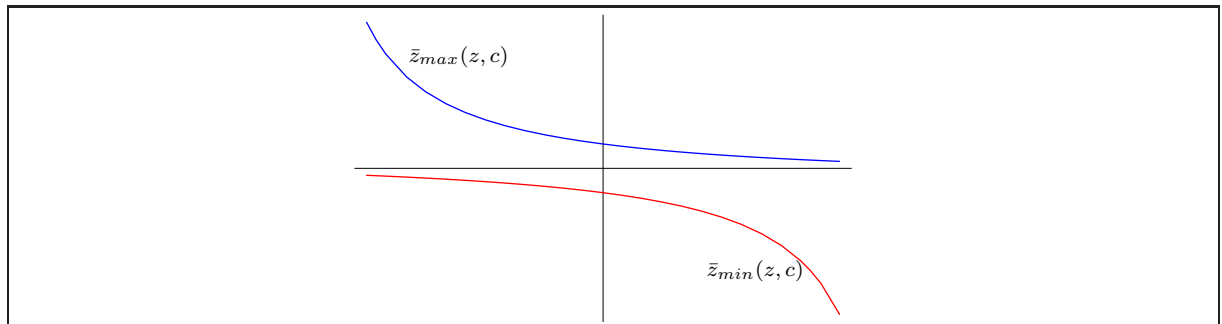}
\caption{The functions $\zbar_{min}(z,c)$ and $\zbar_{max}(z,c)$}\label{F:zbarmin}
\end{figure}

Notice that the set 
\[
\Sigma^c = \scz - (\scz\cap\ec)
\]
is invariant with respect to $Q_z$, $\sigma_{xz}$, $\sigma_{yz}$ and
transpositions of $x$ and $y$. Clearly, $Q_x\Sigma^c$ and
$Q_y\Sigma^c$ do not intersect $\ec$. Therefore
$\Gamma\cdot\Sigma^c \cap \ec= \varnothing$.  Moreover by definition
$|z(u)|>2$ for every $u\in\Sigma^c$ and hence the set $\Lambda_{x,y}$
lies on a hyperbola $\hcz$. Recall that the action of $\Lambda_{x,y}$
on $\hcz$ is wandering and therefore $\Lambda_{x,y}u$ is nowhere
dense in $\hcz$. Consequently, the set
\[
\coprod_{\gamma\in\Gamma}\gamma\Sigma^c
\]
has measure $0$ in $\kappa\inv(c)$. We now have a complete picture
of the $\Gamma$-action on characters for the case when $c<-14$.
\begin{thm}\label{T:ergod}
Suppose $c<-14$. Then the action of $\Gamma$ on the level sets
$\kappa\inv(c)$ is
\begin{enumerate}
\item wandering on the set $\Omega^M$ of characters of discrete
$M$-embeddings 
\item ergodic on the complement of $\Omega^M\cap \kappa\inv(c)$
\[
\kappa\inv(c)\cap(\cv{1}{1} - \Omega^M)
\]
\end{enumerate}
\end{thm}
Next, we extend this argument to the case when
$-14< c < 2$. In that case the inclusion (\ref{E:zminineq}) is
no longer valid. However 
\[
\zbar_{min}(c,z) < 2, \quad \zbar_{max}(c,z)> -2 
\]
for each $z\in(-2,2)$. 
Moreover, there exists a countable open cover
of the interval $(-2, 2)$ by intervals
\[
I_n = (\zbar_{min}(c,z_n), \zbar_{max}(c, z_n))\cap (-2,2),\quad n=1,2,\dots
\]
(see Figure \ref{F:zbarmin}). 
In each $I_n$ a $Q_z$-invariant function $g(z)$ is 
constant almost everywhere by the argument above, and the values of
$g(z)$  agree by default on the overlaps $I_n\cap I_{n+1}$. Therefore
$g(z)$ must be constant almost everywhere in $(-2,2)$. 
\begin{cor}\label{C:ergo}
For any $c<2$ the action of $\Gamma$ on $\kappa\inv(c)\cap(\cv{1}{1} - \Omega^M)$ is
ergodic. 
\end{cor}
We conclude with a remark about the special case $c=-2$. In that case
the level set $\kappa\inv(c)$ consists of characters  $u=(x,y,z)$ that
satisfy the equation
\[
-x^2 - y^2 + z^2 +xyz =0 
\]
The coordinates of such characters are closely related to Markoff
triples, which play an important role in Number Theory. In a recent
paper, Bowditch (\cite{Bowditch}) studies the action of the modular
group on complex Markoff triples and proves that the action has dense
orbits in a neigborhood of the origin. Thus Bowditch's result,
restricted to real characters, is a special case of Corollary~\ref{C:ergo}.

\nocite{*}
\bibliographystyle{amsplain}
\bibliography{paper}

\providecommand{\bysame}{\leavevmode\hbox to3em{\hrulefill}\thinspace}
\begin{thebibliography}{10}

\bibitem{Bowditch}
B.~H. Bowditch, \emph{Markoff triples and quasi-fuchsian groups}, Proc. London
  Math. Soc. \textbf{77} (1998), no.~3, 697--736.

\bibitem{Buser}
P.~Buser, \emph{Geometry and spectra of compact riemann surfaces}, Progress in
  Mathematics, no. 106, Birkh\"auser Boston, 1992.

\bibitem{CullerShalen}
M.~Culler and P.~Shalen, \emph{Varieties of group representations and
  splittings of 3-manifolds}, Ann.\ Math. (1983), no.~117, 109--146.

\bibitem{FrickeKlein}
R.~Fricke and F.~Klein, \emph{Vorlesungen der automorphen funktionen}, vol. I
  (1897) and II (1912), Teubner, Leipzig.

\bibitem{TopComp}
W.~Goldman, \emph{Topological components of spaces of representations}, Inv.
  Math. \textbf{3} (1988), no.~93, 557--607.

\bibitem{Erg}
\bysame, \emph{Ergodic theory on moduli spaces}, Ann. Math. (1997), no.~146,
  475--507.

\bibitem{Expo}
\bysame, \emph{An exposition of results of {F}ricke}, unpublished, May 2000.

\bibitem{Gold:pTorus}
\bysame, \emph{The modular group action on real $\o{SL}(2)$-characters of a
  punctured torus}, preprint, 2002.

\bibitem{Horowitz}
R.~Horowitz, \emph{Induced automorphisms on {F}ricke characters of free
  groups}, Trans.\ A.M.S. (1975), no.~208, 41--50.

\bibitem{}
N.~Ivanov, \emph{Mapping class groups}, Handbook of Geometric Topology (R.J.
  Daverman and R.B. Sher, eds.), Elsevier, 2002, pp.~523--634.

\bibitem{Keen}
L.~Keen, \emph{On {F}ricke moduli}, Advances in the theory of Riemann surfaces
  (Proc. Conf., Stony Brook, N.Y., 1969), Ann. of Math. Studies, no.~66,
  Princeton Univ. Press, 1971, pp.~205--224.

\bibitem{Magnus}
W.~Magnus, \emph{Rings of {F}ricke characters and automorphism groups of free
  groups}, Math.\ Zeit. (1980), no.~170, 91--103.

\bibitem{Magnus-KS}
W.~Magnus, A.~Karrass, and D.~Solitar, \emph{Combinatorial group theory:
  Presentations of groups in terms of generators and relations}, Dover
  Publications, New York, 1977.

\bibitem{Moise}
E.~Moise, \emph{Geometric topology in dimensions $2$ and $3$}, Springer-Verlag,
  Berlin-Heidelberg-New York, 1977.

\bibitem{Nielsen}
J.~Nielsen, \emph{Die isomorphismen der allgemeinen unendlichen gruppe mit zwei
  erzeugenden}, Math.\ Ann. (1918), no.~71, 385--397.

\bibitem{Nielsen1927}
\bysame, \emph{Untersuchungen zur topologie der geschlossenen zweiseitigen
  fl\"achen i}, Acta Math. (1927), no.~50, 189--358.

\bibitem{Procesi}
C.~Procesi, \emph{Invariant theory of $n$ by $n$ matrices}, Adv.\ Math. (1976),
  no.~19, 306--381.

\bibitem{George}
G.~Stantchev, \emph{Action of the modular group on $\glr$-characters of the
  once-punctured torus}, Ph.D. thesis, University of Maryland, 2003.

\bibitem{Steen}
R.~Steenrod, \emph{The topology of fibre bundles}, Princeton Univ. Press, 1951.

\bibitem{Stillwell}
J.~Stillwell, \emph{Classical topology and combinatorial group theory},
  Graduate Texts in Mathematics, Springer-Verlag, New York Berlin Heidelberg,
  1993.

\bibitem{Thurston}
W.~Thurston, \emph{Three-dimensional geometry and topology}, Princeton Univ.
  Press, 1997.

\bibitem{Xia}
E.~Xia, \emph{The moduli of flat $\pglr$ structures on riemann surfaces}, Ph.D.
  thesis, University of Maryland at College Park, 1997.

\end{thebibliography}
\end{document}